\documentclass{elsart}

\usepackage[latin1]{inputenc}
\usepackage{amssymb}
\usepackage{amsfonts}
\usepackage{amsmath}
\usepackage{latexsym}
\usepackage{graphicx,psfrag,epsfig}
\usepackage[numbers]{natbib}
\usepackage[english]{babel}

\newtheorem{theorem}{Theorem}

\newtheorem{corollary}{Corollary}

\newtheorem{lemma}{Lemma}
\newtheorem{proposition}{Proposition}
\newtheorem{remark}{Remark}
\newenvironment{proof}[1][Proof]{\textbf{#1.} }{\ \rule{0.5em}{0.5em}}

\newenvironment{resume}{\noindent{\bf R\'esum\'e}\\
\
\\
\noindent}

\newcommand{\pr}{\mathbb{P}}
\newcommand{\esp}{\mathbb{E}}
\newcommand{\var}{\mathbb{V}\text{ar}}
\newcommand{\sumsum}{\mathop{\sum\sum}}

\newcommand{\betasup}{\bar\beta}
\newcommand{\betainf}{\underline {\beta}}
\newcommand{\rsup}{\bar r}
\newcommand{\rinf}{\underline {r}}
\newcommand{\smin}{\underline s}
\newcommand{\smax}{\bar s}
\newcommand{\alphamin}{\underline \alpha}
\newcommand{\alphamax}{\overline \alpha}
\newcommand{\sgrid}{s_n(s)} 

\def\1{1\!{\mathrm l}}

\newcommand{\Lset}{\mathbb{L}}

\begin{document}

\begin{frontmatter}


\title{Adaptive procedures in convolution models with known or partially known
noise distribution}

\author[l1,l2]{Cristina Butucea}
\address[l1]{\it Laboratoire de Probabilit\'es et
Mod\`eles  Al\'eatoires (UMR CNRS 7599), Universit\'e
Paris VI, 4, pl.Jussieu, Bo\^{\i}te courrier 188, 75252 Paris, France}
\address[l2]{\it Modal'X, Universit\'e Paris X,  200,
avenue de la R\'epublique 92001 Nanterre Cedex, France}
\ead{butucea@ccr.jussieu.fr} 
\author[l3]{Catherine Matias\corauthref{cor1}}
\address[l3]{\it Laboratoire Statistique  et G\'enome (UMR CNRS 8071), Tour
  Evry  2, 523 pl.   des Terrasses  de l'Agora,  91~000 Evry,  France}
\ead{matias@genopole.cnrs.fr} 
\corauth[cor1]{Corresponding author}
\author[l4]{Christophe Pouet}
\address[l4]{\it Laboratoire d'Analyse, Topologie, Probabilit\'es
(UMR CNRS  6632), Centre de  Math\'ematiques et Informatique,  Universit\'e de
Provence, 39 rue F. Joliot-Curie, 13453 Marseille cedex 13, France} 
\ead{pouet@cmi.univ-mrs.fr}

\begin{abstract}
 In a convolution model, we observe random variables whose distribution is the 
 convolution of some unknown density $f$ and some known or partially known 
 noise density $g$.  In this paper, we focus on statistical 
   procedures, which 
 are adaptive with respect to the smoothness parameter $\tau$ of unknown 
 density $f$, 
 and also (in some cases) to some unknown parameter of the noise density $g$.  

In a first part, we assume that $g$ is known and polynomially smooth.  We 
provide goodness-of-fit procedures for the test 
$H_0:f=f_0$, where the alternative $H_1$ is expressed with respect to 
$\Lset_2$-norm (i.e.  has the form $\psi_{n}^{-2}\|f-f_0\|_2^2 \ge 
\mathcal{C}$).  Our adaptive (w.r.t $\tau$) procedure behaves differently 
according to whether $f_0$ is polynomially or exponentially 
smooth.  A payment for adaptation is noted in both cases and for computing 
this, we provide a non-uniform Berry-Esseen type theorem for degenerate 
$U$-statistics. In the first case we prove that the payment for adaptation is 
optimal (thus unavoidable).  

In a second part, we study a wider framework: a semiparametric 
  model, where $g$ is exponentially smooth and stable, 
and its self-similarity index $s$ is unknown.  
 In order to ensure identifiability, we restrict our 
 attention to polynomially smooth, Sobolev-type densities $f$. 
 In this context, we provide 
 a consistent estimation procedure for $s$.  This estimator is then 
 plugged-into three different procedures: estimation of 
 the unknown density $f$, of the functional $\int f^2$ 
 and test of the hypothesis $H_0$. These procedures are 
 adaptive with respect to both $s$ and $\tau$ and attain the rates which are 
 known optimal for known values of $s$ and $\tau$. 
 As a by-product, when the noise is known and exponentially 
 smooth our testing procedure is adaptive for testing Sobolev-type 
 densities.  
\\
\

\begin{resume}
Dans  un  mod\`ele  de  convolution,  les  observations  sont  des  variables
al\'eatoires r\'eelles dont la distribution est la convolu\'ee entre une
densit\'e inconnue $f$ et une densit\'e de bruit $g$ suppos\'ee soit enti\`erement
connue,  soit   connue  seulement  \`a  param\`etre   pr\`es.  Nous  \'etudions
diff\'erentes  proc\'edures  statistiques  qui s'adaptent  automatiquement  au
param\`etre de r\'egularit\'e $\tau$ de la densit\'e inconnue $f$ ainsi que (dans
certains cas), au param\`etre inconnu de la densit\'e du bruit. 

Dans  une  premi\`ere  partie,  nous  supposons  que  $g$  est  connue  et  de
r\'egularit\'e polynomiale. Nous proposons un test
d'ad\'equation de l'hypoth\`ese $H_0: f=f_0$ lorsque l'alternative $H_1$ est
exprim\'ee \`a partir de la norme $\Lset_2$ (i.e. de la forme $\psi_{n}^{-2}\|f-f_0\|_2^2   \ge
\mathcal{C}$).  Cette  proc\'edure est adaptative (par rapport  \`a $\tau$) et
pr\'esente diff\'erentes 
vitesses de test ($\psi_n$) en fonction du type de r\'egularit\'e de $f_0$
(polynomiale ou bien exponentielle).  L'adaptativit\'e induit une perte sur la
vitesse de test, perte qui est calcul\'ee gr\^ace \`a un th\'eor\`eme de type
Berry-Esseen non-uniforme  pour des $U$-statistiques  d\'eg\'en\'er\'ees. Dans
le  cas d'une  r\'egularit\'e polynomiale  pour $f$,  nous prouvons  que cette
perte est in\'evitable et donc optimale. 

Dans un second temps, nous nous placons dans le cadre plus large d'un mod\`ele
semi-param\'etrique,   o\`u   $g$   est   la  densit\'e   d'une   loi   stable
(r\'egularit\'e de  type exponentiel)  avec un indice  d'auto-similarit\'e $s$
inconnu. Pour  assurer l'identifiabilit\'e du  mod\`ele, la densit\'e  $f$ est
suppos\'ee   appartenir    \`a   un   espace    de   Sobolev   (r\'egularit\'e
polynomiale).   Dans ce  cadre,  nous proposons  un  estimateur consistant  de
$s$.  Celui-ci est ensuite inject\'e dans trois proc\'edures diff\'erentes :
l'estimation  de  $f$,   de  la  fonctionnelle  $\int  f^2$   et  le  test  de
l'hypoth\`ese  $H_0$. Ces proc\'edures  sont adaptatives  par rapport  \`a $s$
et  \`a $\tau$  et atteignent  les  vitesses optimales  du cas  $s$ et  $\tau$
connus.  Enfin, lorsque $g$ est connue et de r\'egularit\'e
exponentielle, une cons\'equence de notre r\'esultat est que cette proc\'edure
de test est adaptative lorsque $f_0$ appartient \`a un espace de Sobolev. 
\end{resume}
\end{abstract}

\end{frontmatter}


\section{Introduction}
\subsection*{Convolution model}
Consider the \textbf{convolution model} where the observed sample $\{Y_j\}_{1\leq j
\leq  n}$  comes from  the  independent  sum  of independent  and  identically
distributed  (i.i.d.) random  variables  $X_j$ with  unknown  density $f$  and
Fourier  transform  $\Phi$ and  i.i.d.  noise  variables $\varepsilon_j$  with
known (maybe only up to a parameter) density $g$ and Fourier transform $\Phi^g$
\begin{equation}\label{model}
 Y_j = X_j + \varepsilon_j, \quad 1\leq j \leq n.
\end{equation}
The density of the observations is denoted by $p$ and its Fourier transform
$\Phi^p$. Note that we have $p = f *g$ where $*$ denotes
the convolution product and $\Phi^p = \Phi \Phi^g$.

The underlying unknown density $f$ is always supposed to belong to $\mathbb{L}_1 \cap
\mathbb{L}_2$. We shall consider probability density functions belonging to the class
\begin{equation}
\mathcal{F} \left( \alpha, r, \beta, L\right) =\left\{ f:\mathbb{R}\rightarrow \mathbb{R}_{+}
, \int f =1 , \frac{1}{2 \pi}\int \left| \Phi \left(
u\right) \right| ^{2} |u|^{2\beta} \exp \left( 2\alpha |u|^{r}\right) du\leq L\right\},
\label{mixte}
\end{equation}
for $ L$ a positive constant, $\alpha > 0$, $0 \leq r \leq 2$, $\beta \geq 0$ and
either $r>0$  or $\beta >0$. Note  that the case $r=0$  corresponds to Sobolev
densities whereas $r>0$ corresponds to infinitely many differentiable (or supersmooth) densities.

We consider noise distributions whose Fourier transform  does not vanish
on $\mathbb{R}$:  $\Phi^g(u) \neq 0$, $\forall~u  \in \mathbb{R}$.  Typically,
nonparametric estimation in
convolution models gives rise to  the distinction of two different behaviours
for the noise distribution. We alternatively shall consider (for some constant
$c_g>0$),

\textbf{polynomially smooth} (or polynomial) noise
\begin{equation}
\left| \Phi ^{g}\left( u\right) \right| \sim c_g \left| u\right| ^{-\sigma}\text{, }%
\left| u\right| \rightarrow \infty \text{, } \sigma >1\text{;}
\label{polynomial}
\end{equation}

\textbf{exponentially  smooth  }  (or  supersmooth  or
exponential) \textbf{stable} noise 
\begin{equation}
\left| \Phi ^{g}\left( u\right) \right| = \exp \left( -\gamma \left|
u\right| ^{s}\right) \text{, }\left| u\right| \rightarrow \infty \text{, }%
\gamma ,s>0\text{.}
\label{exponential}
 \end{equation}
In this second case, the parameter $s$ is called the
self-similarity  index of  the noise  density and  we shall consider that it  is {\bf
  unknown}.

Convolution models have been widely studied over the past two decades.
We will be interested here both in estimation of the unknown density $f$ and
in testing the hypothesis $H_0: f=f_0$,  with  a  particular  interest  in  {\it  adaptive}
procedures.  Our first purpose is to provide goodness-of-fit
testing  procedures on $f$,  for the  test of  the hypothesis
  $H_0: f=f_0$, which 
are adaptive with respect to the unknown smoothness
parameter of $f$.  The second one is to  study the behaviour
of different procedures (such 
as estimation of $f$, estimation of  $\int f^2$ and goodness-of-fit test) in a
setup where self-similarity index $s$ is unknown.

\subsection*{Adaptive procedures in the convolution model}
Concerning estimation, the asymptotically minimax setup in the
context  of pointwise or  $\Lset_p$-norms and  in the  case of
entirely known noise density $g$ is the most studied one. Major
results in this direction prove that the smoother the error density,
the slower       the      optimal       rates   of convergence
(see \citep{Carroll-Hall}, \citep{Fan1}, \citep{Butucea-deconv},
\citep{Fan3} concerning
polynomial noise and \citep{penvid}, \citep{Comte-Taupin-Rozen},
\citep{ButTsy} for exponential noise).
Adaptive estimation  procedures were considered first by
\citep{penvid} and then by \citep{Fan-Koo}. They constructed wavelets
estimators which do not depend on smoothness parameter of the
density   $f$ to be estimated. Adaptive kernel estimators were given
in \citep{ButTsy}. A different adaptive approach is used in
\citep{Comte-Taupin-Rozen} relying on penalized contrast estimators.

Nonparametric goodness-of-fit testing has extensively been studied in the context of direct
observations (namely a sample distributed  from the density $f$ to be tested),
but also  for regression  or in the  Gaussian white  noise model. We  refer to
\citep{LehRom}, \citep{Ingster} for an overview on the subject.  The
convolution model  provides an interesting  setup where observations  may come
from a signal observed through some noise.
\\

Nonparametric goodness-of-fit tests in convolution models were
studied in \citep{HolzBisMunk} and in \citep{Butucea}. The
approach used in \citep{Butucea} is based on a minimax point of
view combined with estimation of the quadratic functional $\int
f^2$. Assuming the  smoothness parameter  of $f$  to be  known,
the authors of \citep{HolzBisMunk} define a version of the Bickel-Rosenblatt
test statistic and study its
asymptotic distribution under the null hypothesis and under fixed
and local alternatives, while \citep{Butucea} provides a
different goodness-of-fit testing procedure  attaining the  minimax
rate  of testing  in each  of  the three following  setups: Sobolev
densities  and   polynomial  noise,  supersmooth densities and
polynomial noise,  Sobolev densities and exponential noise. The case
of supersmooth densities and  exponential noise is  also studied but
the optimality of the procedure is not established in the case $r>s$.\\

 Our first goal here is to provide adaptive versions of these last procedures with
respect to  the parameters $(\alpha,r,\beta)$.   We restrict our  attention to
testing   problems  where   alternatives   are  expressed   with  respect   to
$\Lset_2$-norm.   Namely, the alternative has the form $H_1:
\psi_{n}^{-2} \|f-f_0\|_2^2 \ge
\mathcal{C}$. In such a case, the problem relates with
asymptotically minimax estimation of $\int f^2$.  \\

 Our second goal is to deal with the case of not entirely known noise
distribution. This is  a crucial issue as, assuming  this noise distribution
  to be entirely known is not realistic in many situations. However, 
in  general, the  noise  density $g$  has  to be  known for  the
model to  be identifiable. Nevertheless, when the noise density is
exponentially smooth and the  unknown density is  restricted to be
less smooth than  the noise, semiparametric  models are identifiable
and they may be considered. The  case of a
Gaussian noise  with unknown  variance $\gamma$  and unknown density
without   Gaussian   component   has  first   been   considered   in
\citep{Matias}. She  proposes an estimator  of the parameter $\gamma$
which is then plugged in an estimator of  the unknown density. This
work is generalized in \citep{Butmat} for exponentially smooth noise
with unknown scale parameter $\gamma$ and unknown  densities
belonging either to Sobolev  classes, or to classes of supersmooth densities
with parameter $r$, $r<s$.  Minimax rates of convergence are exhibited.
In  this context, the unknown parameter  $\gamma$ acts as  a real
nuisance parameter  as the  rates of convergence  for estimating the unknown
density are slower compared  to the  case of known  scale, those
rates  being nonetheless optimal in a minimax sense.  Another
attempt to remove knowledge on the noise density appears in
\citep{Meister} where the author studies a deconvolution estimator
associated to a procedure for selecting the error density
between the Normal supersmooth density and the Laplace polynomially smooth
density (both with fixed parameter values). \\

 In the second part of our work, we will be interested in estimation procedures on $f$,
adaptive  both with  respect to  the  smoothness parameter  of $f$  and to  an
unknown parameter of the noise density. More precisely, in the specific
setup of Sobolev densities and exponential noise with symmetric stable distribution,
we will consider the case of
unknown self-similarity index $s$. In this context, we first propose an estimator of
the self-similarity index $s$, which, plugged into kernel procedures, provides
estimators of the unknown density  $f$ with the same optimal rate  of convergence as in
the case of  entirely known noise density. Using the  same techniques, we also
construct an  estimator of the  quadratic functional $\int f^2$  (with optimal
rate of convergence) and $\mathbb{L}_2$ goodness-of-fit test statistic.
Note that this  work is very different from \citep{Butmat} as the self
similarity index $s$ plays a  different role from the scale parameter $\gamma$
previously studied. In particular, the  range of applications of those results
is entirely new.

\subsection*{Notation, definitions, assumptions}

In  the  sequel,  $\|\cdot\|_2$  denotes  the $\mathbb{L}_2$-norm,
$\bar  M$ is the  complex conjugate of  $M$ and $<M  ,N > =\int
M(x)\bar N(x) dx$ is the scalar product of complex-valued functions
in $\mathbb{L}_2(\mathbb{R})$. Moreover, probability and expectation
with respect to the distribution of $Y_1, \ldots ,Y_n$ induced by
the unknown density $f$ will
be denoted by $\pr_f$ and $\esp_f$.\\

We denote more generally  by $\tau=(\alpha,r, \beta)$ the smoothness parameter
of the unknown density $f$ and by $\mathcal{F}(\tau,L)$ the corresponding class.
As the density $f$ is unknown, the a priori knowledge of its smoothness
parameter $\tau$ could appear unrealistic.
Thus, we assume that $\tau$ belongs to a closed subset $\mathcal{T}$, included
in $(0,+\infty) \times (0,2] \times (0,+\infty)$.
For a  given density $f_0$  in the class  $ \mathcal{F}(\tau_0)$, we  want to
test the hypothesis
$$ H_0 : f = f_0 $$
from observations $Y_1,\ldots,Y_n$ given by \eqref{model}.
We  extend the  results of  \citep{Butucea} by  giving the family of
sequences $\Psi_n = \{\psi_{n,\tau}\}_{\tau  \in \mathcal{T}}$  which separates  (with  respect to
$\mathbb{L}_2$-norm) the null hypothesis from a larger alternative
\begin{equation*}
H_1(\mathcal{C},  \Psi_n)  :  f  \in \cup_{\tau  \in  \mathcal{T}}\{f  \in
\mathcal{F}(\tau,L) \text{
  and } \psi_{n,\tau}^{-2} \|f-f_0\|_2^2 \geq \mathcal{C} \}.
\end{equation*}
We recall that the usual procedure is to construct, for any $0<\epsilon<1$, a test
statistic $\Delta_n^\star$  (an arbitrary function, with  values in $\{0,1\}$,
which is measurable with respect to  $Y_1,\ldots ,Y_n$ and such that we accept
$H_0$ if $\Delta_n^\star=0$ and reject it otherwise) for which there exists some $\mathcal{C}^0 >0$ such that
\begin{equation}\label{uptest}
\limsup_{n \to \infty} \left\{ \pr_{0}[\Delta_n^\star = 1] +
\sup_{f \in H_1(\mathcal{C},\Psi_{n})} \pr_{f}[\Delta_n^\star =
0]\right\} \leq \epsilon,
\end{equation}
holds for all $\mathcal{C}>\mathcal{C}^0$.  This part is called the upper
bound of the testing rate.
Then, prove the minimax optimality of this procedure, i.e. the lower bound
\begin{equation} \label{lowtest}
\liminf_{n \to \infty} \inf_{\Delta_n} \left\{ \pr_{0}[\Delta_n = 1]
+\sup_{f \in H_1(\mathcal{C}, \Psi_{n})} \pr_{f}[\Delta_n = 0]\right\} \geq \epsilon,
\end{equation}
for  some  $\mathcal{C}_0>0$ and for all  $0<\mathcal{C}<\mathcal{C}_0$, where the
infimum is taken over all test statistics $\Delta_n$.

Let us first remark that as we use noisy observations (and unlike what happens
with direct observations), this test cannot be
reduced  to testing  uniformity of  the distribution  density of  the observed
sample (i.e. $f_0 =1$ with support on the finite interval $[0;1]$).
As a consequence, additional assumptions used in \citep{Butucea} on the
tail  behaviour of $f_0$  (ensuring it  does not  vanish arbitrarily  fast) are
needed to obtain the optimality result of the testing procedure in the case of
Sobolev density ($r=0$) observed with polynomial noise ({\bf (T)} and {\bf (P)}),
respectively with exponential noise ({\bf (T)} and {\bf (E)}).
We recall these assumptions here for reader's convenience.

{\bf Assumption (T)}
  $$ \exists c_0>0, \forall x \in \Rset , f_0(x) \geq \frac{c_0}{1+|x|^2} .$$

Moreover, we also need to control the derivatives of known Fourier transform
$\Phi^g$ when establishing optimality results.

 {\bf Assumption (P)} (Polynomial noise) If the noise satisfies (\ref{polynomial}),
 then assume that $\Phi^g$ is three times continuously
 differentiable and there exist $A_1, A_2$ such that
 $$ |(\Phi^g )'  (u ) | \leq \frac{A_1}{|u|^{\sigma+1}} \text{  and } |(\Phi^g )^{''}
 (u ) | \leq \frac{A_2}{|u|^{\sigma+2}}, \quad |u| \to \infty .$$

    {\bf Assumption (E)} (Exponential noise) If the noise satisfies (\ref{exponential}),
    then assume that $\Phi^g$ is continuously differentiable and
    there exists some constants $C>0$ and $A_3 \in \Rset$ such that
$$ | (\Phi^g )' (u ) | \leq C |u|^{A_3} \exp(-\gamma |u| ^s) , \quad |u| \to \infty .$$

\begin{remark} Similar results may be obtained when we assume the existence of
  some   $p\geq 1$  such   that  $f_0(x)$  is bounded from below by
  $c_0  (1 +|x|^{p})^{-2}$ for large enough $x$.
  In such  a case, the Fourier transform
  $\Phi^g$ of the  noise density is assumed to be $p$ times continuously differentiable, with
  derivatives up to order $p$ satisfying  the same  kind  of  bounds  as in  Assumption  {\bf
    (P)}, when the noise is polynomial, respectively in Assumption~{\bf (E)},
    when the noise is exponential.
\end{remark}

\subsection*{Roadmap}
Section~\ref{sec:bruit_poly} deals with the  case of (known) polynomial noise.
We provide  a goodness-of-fit testing procedure for the test
  $H_0: f=f_0$, in two different cases: the
density $f_0 $ to be tested is either ordinary smooth ($r_0=0$)
or supersmooth ($r_0>0$). The procedures are adaptive with
respect to the smoothness parameter $(\alpha, r, \beta)$ of $f$.
The proof of the upper bounds for the testing rate relies mainly on a Berry-Esseen inequality
for degenerate $U$-statistics of order $2$, postponed to
Section~\ref{sec:Ustats}. In some cases, a loss for adaptation is noted with respect
to known testing rates for fixed known parameters. When the loss is of order
$\log \log n$ to some power, we prove that this payment is unavoidable.

In  Section~\ref{sec:s_inconnu}, we  consider exponential  noise of symmetric
stable law  with unknown self-similarity index $s$.
In  order to  ensure identifiability,  we  restrict our
attention to Sobolev classes of densities $f$.
The first step (Section~\ref{sec:estim_s}) is
to provide  a consistent estimation  procedure for the  self-similarity index.
Then (Section~\ref{sec:apres_estim_s}) using a plug-in, we introduce a new
kernel estimator of $f$ where both the bandwidth and the kernel are data dependent.
We also introduce an estimator of the quadratic functional $\int f^2$ with
sample dependent bandwidth and kernel. We prove that these two
procedures attain  the same rates  of convergence as  in the case  of entirely
known noise distribution, and are thus asymptotically optimal in the minimax sense.
We also present a goodness-of-fit test on $f$ in this setup.
We prove that the  testing rate is the
same  as  in   the  case  of  entirely  known   noise  distribution  and  thus
asymptotically optimal in the minimax sense.
Proofs are postponed to Section~\ref{sec:proofs}.

%
\section{Polynomially smooth noise}\label{sec:bruit_poly}

In this section, we shall assume that the noise density $g$ is polynomial \eqref{polynomial}.
The unknown density $f$ belongs to the class $\mathcal{F}(\alpha,r,\beta,L)$. 
We are interested in adaptive, with respect to the parameter $\tau =(\alpha,r,
\beta)$,  goodness-of-fit  testing procedures.  We  assume  that this  unknown
parameter belongs to the following set
$$
\mathcal{T} = \{ \tau=(\alpha,r, \beta) ;  \tau \in
[\alphamin;+\infty) \times [\rinf;\rsup]\times [\betainf;\betasup]
\},
$$
where $\alphamin>0$, $0\leq \rinf \leq \rsup \leq 2$, $0\leq
\betainf \leq \betasup$ and either $\rinf >0$ and $\alpha \in [\alphamin,\alphamax]$
or both $\rinf=\rsup=0$ and $\betainf >0$.

\bigskip

Let us introduce some notation.
We consider a preliminary kernel $J$, with Fourier transform
$\Phi^J$, defined by
$$
\forall x \in \mathbb{R}, \; \; J(x) = \frac{\sin(x)}{\pi x} , \quad \forall u
\in \mathbb{R}, \; \; \Phi^J(u) = 1_{|u| \leq 1},
$$
where $1_{A}$ is the indicator function of the set $A$.
For any bandwidth $h= h_n \to 0$ as $n$ tends to infinity, we define
the rescaled kernel $J_h$  by
$$ \forall x \in \mathbb{R}, \; \;
J_{h}(x)  =  h^{-1} J(x/h)  \text{  and  } \forall  u  \in  \mathbb{R}, \;  \;
\Phi^{J_h} (u) = \Phi^{J} (hu) = 1_{|u| \leq 1/h}.
$$
Now, the  deconvolution kernel $K_{h}$ with  bandwidth $h$ is  defined via its
Fourier transform $\Phi^{K_{h}}$ as
\begin{equation}
\Phi^{K_{h}} (u) = \left( \Phi^g (u)\right)^{-1} \Phi^J (uh)
=  \left( \Phi^g (u)\right)^{-1} \Phi^{J_h} (u)  , \quad \forall u \in \mathbb{R}.
\end{equation}

In Section~\ref{sec:apres_estim_s}, we will consider a modification of this
kernel to take into account the case of not entirely known noise density $g$.

Next, the quadratic functional $\int (f-f_0)^2$ is estimated by the statistic $T_{n,h}$
\begin{eqnarray}
T_{n, h} &=& \frac{2}{n(n-1)} \sumsum_{1\leq k < j \leq n} <K_{h}(\cdot -Y_k) -f_0, K_{h}(
\cdot -Y_j) -f_0 >.  \label{Tstat}
\end{eqnarray}
Note that $T_{n,h}$ may not be positive, but its expected value is.

\bigskip

In order to construct a testing procedure which is adaptive with respect to the parameter $\tau$
we introduce a sequence of finite regular grids over the set $\mathcal{T}$ of
unknown parameters: $\mathcal{T}_N =\{\tau_i ; 1\leq i \leq N \} $.
For each grid point $\tau_i$ we choose a testing threshold $ t_{n,i}^2$ and a
bandwidth $h_n^i$ giving a test statistic $T_{n,h_n^i}$.

The test rejects the null hypothesis  as soon as at least
one  of  the  single tests  based  on  the parameter  $\tau_i$  is
rejected.
\begin{equation}\label{test}
\Delta_n^\star =\left\{
\begin{array}{ll}
1 & \text{if } \sup_{1\leq i\leq N} | T_{n,h_n^i} | t_{n,i}^{-2}
>\mathcal{C}^\star\\
0 & \text{otherwise},
\end{array}
\right.
\end{equation}
for  some constant  $\mathcal{C}^\star>0$ and  finite sequences  of bandwidths
$\{h_n^i  \}_{1\leq i  \leq N}$  and thresholds  $\{t_{n,i}^2\}_{1\leq  i \leq
  N}$.

We  note  that   our  asymptotic  results  work  for   large  enough  constant
$\mathcal{C}^\star$.
In practice we may choose it by Monte-Carlo simulation under the null hypothesis,
for known $f_0$, such that we control the first-type error of the test
and bound it from above, e.g. by $\epsilon/2$.

Typically, the structure of the grid accounts for two different
phenomena. A first part of the points is dedicated to the adaptation
with respect to  $\beta$ in  case $\rsup=\rinf=0$, whereas the rest
of the points is used to adapt the procedure with respect to $r$
(whatever the value of $\beta$).

In the two next theorems, we fix $\sigma >1$. We note that the testing
rates are essentially different according to the two different cases where $f_0$
belongs to a Sobolev class ($r_0 =0$, $\alpha_0 \geq \alphamin$
and we assume $\beta_0=\betasup$) and where $f_0$ is a supersmooth
function  ($\alpha_0 \in [\alphamin,\alphamax]$, $r_0>0$ and $\beta_0 \in [\betainf,\betasup]$
and then  we focus on $r_0=\rsup$ and $\alpha_0=\alphamax$).  Note that in the first case, the
alternative contains functions $f$ which are smoother ($r>0$) than
the null hypothesis $f_0$.

\bigskip

When $f_0$  belongs to Sobolev  class $\mathcal{F}(\alpha_0,0,\betasup,L)$, the
grid is defined as follows. Let $N$ and choose $\mathcal{T}_N =\{\tau_i ; 1\leq i \leq N+1 \}$ such that
\begin{multline*}
\left\{
\begin{array}{l}
\forall 1\leq i \leq N,  \tau_i=(0;0;\beta_i) \text{ and }
\beta_1=\betainf  <\beta_2 <\ldots < \beta_{N} =\betasup,\\
\forall  1\leq i \leq N-1, \; \beta_{i+1} -\beta_i = (\betasup-\betainf)/ (N-1) ,  \\
\text{ and } \tau_{N+1} =(\alphamin;\rsup;0)
\end{array}
\right.
\end{multline*}
In this case, the first $N$ points  are  dedicated to  the
adaptation with respect to $\beta$ when $\rsup=\rinf=0$,
whereas the last point $\tau_{N+1}$ is  used to  adapt the procedure
with respect to $r$ (whatever the value of $\beta$).

\begin{theorem}\label{th:bruit_poly_f0sob}
Assume $f_0 \in \mathcal{F}(\alpha_0, 0, \betasup, L)$.
The test statistic $\Delta_n^\star$ given by \eqref{test} with parameters
\begin{eqnarray*}
 N= \lceil \log n \rceil;  \quad
\forall 1\leq i \leq N :
\left\{
\begin{array}{l}
h_n^i= \left( \frac{n}{\sqrt{\log \log n}}\right)^{-2/(4\beta_i +4\sigma+1)}\\
t_{n,i}^2 = \left( \frac{n}{\sqrt{\log \log n}}\right)^{-4\beta_i/(4\beta_i +4\sigma+1)}
\end{array}
\right. , \\
h_n^{N+1}= n^{-2/(4 \betasup +4\sigma+1)} ; \quad t_{n,N+1}^2 =
n^{-4\betasup /(4 \betasup +4\sigma+1)} ,
\end{eqnarray*}
and any large enough positive constant $\mathcal{C}^\star$, satisfies \eqref{uptest} for any
$\epsilon \in (0,1)$, with testing rate $\Psi_n = \{\psi_{n, \tau} \}_{\tau \in \mathcal{T}}$ given by
$$
\psi_{n,\tau} =
 \left( \frac{n}{\sqrt{\log \log n}}\right)^{-2 \beta/(4 \beta +4\sigma+1)} 1_{r=0}
 + n^{-2 \betasup/(4 \betasup +4\sigma+1)} 1_{r>0},\, \forall\, \tau=(\alpha,r,\beta) \in \mathcal{T}.
$$
Moreover, if $f_0 \in \mathcal{F}(\alpha_0, 0, \betasup, cL)$ for
some $0<c<1$ and if Assumptions
{\bf (T)} and  {\bf(P)} hold, then this testing rate  is adaptive
minimax over the family of classes $\{ \mathcal{F}(\tau, L), \tau
\in [\alphamin,\infty) \times \{0\}\times [\betainf,\betasup]\}$
(i.e. \eqref{lowtest} holds).
\end{theorem}

We note that our testing procedure attains the polynomial rate
$n^{-2 \betasup/(4 \betasup +4\sigma+1)}$ over the union of all classes
containing functions smoother than $f_0$.
Note moreover that this rate is known to be a minimax testing rate
over the class $\mathcal{F}(0,0,\betasup,L)$ by results in
\citep{Butucea}. Therefore we prove that the loss of some power of
$\log \log n$ with respect to the minimax rate is unavoidable. A
loss appears when the alternative contains classes of functions less
smooth than $f_0$.

The proof that  our adaptive procedure attains the minimax  rate relies on the
Berry-Esseen inequality presented in Section~\ref{sec:Ustats}.

\bigskip

When $f_0$ belongs to class $\mathcal{F}(\alphamax,\rsup,\beta_0,L)$ of infinitely many
differentiable functions, the  grid is defined as follows.  Let $N_1, \, N_2$
and choose
$\mathcal{T}_N =\{\tau_i ; 1\leq i \leq N=N_1+N_2 \}$ such that
\begin{multline*}
\left\{
\begin{array}{l}
\forall 1\leq i \leq N_1,  \tau_i=(0;0;\beta_i) \text{ and }
\beta_1=\betainf  <\beta_2 <\ldots < \beta_{N_1} =\betasup,\\
\forall  1\leq i \leq N_1-1, \; \beta_{i+1} -\beta_i = (\betasup-\betainf)/ (N_1-1) ,  \\
\text{ and } \forall 1\leq i \leq N_2,  \tau_{N_1+i}=(\alphamax;
r_i;\beta_0) \text{ and }
r_1=\rinf  <r_2 <\ldots < r_{N_2} =\rsup,\\
\forall  1\leq i \leq N_2-1, \; r_{i+1} -r_i = (\rsup-\rinf)/
(N_2-1) .
\end{array}
\right.
\end{multline*}
In this case, the first $N_1$ points are used for  adaptation with
respect to  $\beta$ in  case $\rsup=\rinf=0$, whereas the last
$N_2$ points are  used to  adapt the procedure with respect to $r$
(whatever the value of $\beta$).

\begin{theorem}\label{th:bruit_poly_f0ana}
Assume  $f_0 \in \mathcal{F}(\alphamax, \rsup, \beta_0, L)$ for some
$\beta_0 \in [\betainf,\betasup]$.
The test statistic $\Delta_n^\star$ given by
\eqref{test} with $\mathcal{C}^\star$ large enough and
$$ N_1= \lceil \log n \rceil ;  \quad
\forall 1\leq i \leq N_1 :  \left\{ 
\begin{array}{l}
h_n^i= \left( \frac{n}{\sqrt{\log \log n}}\right)^{-2/(4\beta_i +4\sigma+1)}\\
t_{n,i}^2 = \left( \frac{n}{\sqrt{\log \log n}}\right)^{-4\beta_i/(4\beta_i +4\sigma+1)}
\end{array}
\right. , 
$$
$$
 N_2=  \lceil \log \log n /(\rsup-\rinf) \rceil  ;  
\forall 1\leq i \leq N_2 :  \left\{ 
\begin{array}{l}
h_n^{N_1+i}= \left( \frac{\log n}{2c}\right)^{-1/r_i } ,
c < \alphamin \exp\left( -\frac{1}{\rinf}\right)\\
t_{n,N_1+i}^2 = \frac{\left(\log n\right)^{(4\sigma+1)/(2 r_i)}}{n}
\sqrt{\log \log \log n}
\end{array} ,
\right. 
$$
satisfies \eqref{uptest}, with testing rate $\Psi_n = \{\psi_{n,
\tau} \}_{\tau \in \mathcal{T}}$ given by
$$
\psi_{n,\tau}     =   \left( \frac{n}{\sqrt{\log \log n}}\right)^{-2
\beta/(4 \beta +4\sigma+1)} 1_{r=0}
 + \frac{\left(\log n \right)^{(4\sigma+1)/(4 r)}}{\sqrt{n}}(\log \log \log n)^{1/4}  1_{r\in [\rinf,\rsup]}  .
$$
\end{theorem}

We note that if Assumptions {\bf (T)} and  {\bf(P)} hold for $f_0$
in $\mathcal{F}(\alphamax, \rsup, \beta_0, L)$, the same optimality
proof as in Theorem~\ref{th:bruit_poly_f0sob} gives us that the loss
of the $\log \log n$ to some power factor is optimal over
alternatives in $\bigcup_{\alpha \in[\alphamin, \alphamax], \beta
\in [\betainf,\betasup]} \mathcal{F}(\alpha, 0, \beta, L)$. A loss
of a $(\log \log \log n)^{1/4}$ factor appears over alternatives of
supersmooth densities (less smooth than $f_0$) with respect to the
minimax rate in \citep{Butucea}. We do not prove that this loss
is optimal.

%
\section{Exponentially smooth noise in a semiparametric context}\label{sec:s_inconnu}

In this section, we assume  the noise  density $g$ to be
exponentially  smooth  and  stable \eqref{exponential},  for
some {\bf unknown} 
$s \in [\smin;\smax]$ and fixed (known) bounds $0<\smin<\smax\leq
2$. More precisely, we suppose that the noise has symmetric stable
law having Fourier transform \\

{\bf Assumption (S)}
$ \; \; \Phi^g (u)  = \exp(-|u|^s) \text{ where  } s \in [\smin; \smax] $. \\

The results of  Section~\ref{sec:estim_s} are valid under the
  more  general  assumption  \eqref{exponential}  with known  scale  parameter
  $\gamma$,  which enables  us to  select the  smoothness parameter  among the
  wider class of  not necessarily symmetric stable densities  with known scale
  parameter. Nevertheless, the exact form of Fourier transform $\Phi^g $ 
is needed for deconvolution purposes (see Section~\ref{sec:apres_estim_s}).

For the model to be identifiable, we must assume that $f$
is not too smooth, i.e. its Fourier transform does not decay  asymptotically faster
than a known polynomial of order $\beta'$. \\

{\bf Assumption (A)}  There exists some known $A>0$,  such that
 $ |\Phi(u)| \geq A|u|^{-\beta'} $  for large enough $|u|$.\\

The notation $q_{\beta'}$ is used for the function $u \mapsto A |u|^{-\beta'}$.
Under assumptions {\bf(S)}  and {\bf(A)} the model is identifiable. Indeed,
considering the Fourier transforms, we get for all real number $u$
$$ \log |\Phi^p(u)| = \log |\Phi(u)| - |u|^{s}.$$
 Now  assume that we
have the equality between two Fourier transform for the observations $\Phi^p_1
= \Phi^p_2$,  where $\Phi^p_1(u) =\Phi_1(u) e^{-|u|^{s_1}}$ and
$\Phi^p_2(u) =\Phi_2(u)  e^{-|u|^{s_2}}$.
Without  loss of generality, we may assume $s_1  \leq s_2$. Then we get
$$
|u|^{-s_1} \log |\Phi_1(u)|-1 = |u|^{-s_1} \log |\Phi_2(u)| -|u|^{s_2-s_1}
$$
and taking  the limit  when $|u|$ tends  to infinity implies  (with assumption
{\bf(A)}) that $s_1=s_2$ and then
$\Phi_1 =\Phi_2$ which proves the identifiability of the model. \\

In this context, $\pr_{f,s}$ and $\esp_{f,s}$ respectively denote probability
and expectation with respect to the model under parameters $(f,s)$.\\

\subsection{Estimation of the self-similarity index $s$}\label{sec:estim_s}
We first present a selection procedure $\hat s_n$ which asymptotically recovers the
true  value of the  smoothness parameter  $s$, with fast rate of convergence.
We use a discrete grid  $\{s_1, \ldots  , s_N\}$, with a number $N$ of points growing to infinity.

\begin{figure}[!htbp]
 \begin{center}
 \psfrag{qphik-1}{$q_{\beta'} \Phi^{[k-1]}$}
 \psfrag{phik}[][]{$\Phi^{[k]}$}
 \psfrag{qphik}[][]{$q_{\beta'} \Phi^{[k]}$}
 \psfrag{phik+1}[][]{$\Phi^{[k+1]}$}
 \psfrag{ICI}[][]
{$ \left\} s_k \phantom{
     \begin{array}{c}
 1 \\ 2 \\3
     \end{array}}
    \right. $}
 \includegraphics[width=12cm]{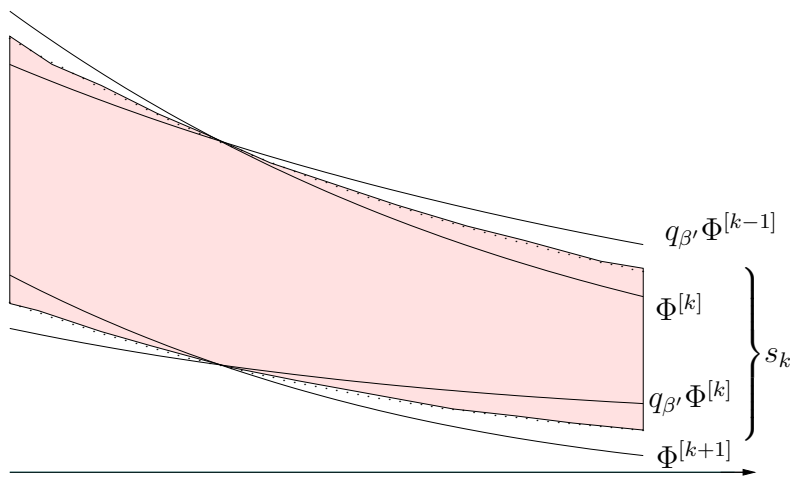}
 \caption{Estimation procedure for $s$.  When $|\hat \Phi^p_n (u_n) |$ lies in
   the grey region, we choose $\hat s_n= s_k$.}
\label{fig}
 \end{center}
\end{figure}

The asymptotic behavior of the  Fourier transform $\Phi^p$ of the observations
is used to  select the smoothness index $s$. More  precisely, we have for
any large enough $|u|$
$$
A|u|^{-\beta'}   \exp(-|u|^s)  \leq  |\Phi^p(u)|   \leq
\exp(-|u|^s),
$$
namely, the function $|\Phi^p|$ asymptotically belongs to the {\it pipe}
$[q_{\beta'}(u)  e^{-|u|^{s}};  e^{-|u|^{s}}]$.
Let us now consider  a discrete grid $0<\smin=s_1<s_2<\ldots <s_N=\smax\leq 2$
and denote $\Phi^{[k]}(u) =  e^{-|u|^{s_k}}$.
These families of functions $\{\Phi^{[k]}\}_{1\leq k \leq N}$ and $\{ q_{\beta'}
\Phi^{[k]}\}_{1\leq  k \leq  N}$ form  an asymptotically  decreasing  family as
there exists some  positive real number $u_1$ such that for  all real $u$ with
$|u| \geq u_1$, we have
\begin{equation}
  \label{order_function}
\Phi^{[1]}(u) \geq  q_{\beta'}(u)\Phi^{[1]}(u) \geq \Phi^{[2]}(u)
\geq  \cdots \geq \Phi^{[N]}(u) \geq q_{\beta'} (u) \Phi^{[N]}(u).
\end{equation}
If the  size of  the grid is  sufficiently small,  the modulus of  the Fourier
transform $\Phi^p$ will asymptotically belong to one of the pipes $[q_{\beta'}
\Phi^{[k]}; \Phi^{[k]}]$.
Our estimation procedure uses the empirical estimator
$$ \hat \Phi^p_n(u) =\frac{1}{n}\sum_{j=1}^n \exp(-iuY_j) , \quad \forall u \in \Rset ,
 $$
of the Fourier  transform $\Phi^p$ at some point $u_n$ which tends to infinity with $n$.
The procedure selects the  smoothness parameter among
$\{s_1, \ldots  , s_N\}$  by choosing the  pipe
$[q_{\beta'}(u_n) \Phi^{[k]}(u_n) ; \Phi^{[k]}(u_n)]$
closest to the function $|\hat \Phi_n^p(u_n)|$. More precisely
\begin{equation}
  \label{estim_s}
  \hat s_n = \left\{
  \begin{array}{ll}
 s_k &\text{if }
\frac{1}{2}  \left\{  q_{\beta'} \Phi^{[k]} +\Phi^{[k+1]}  \right\}  (u_n) \leq  |\hat
\Phi^p_n  (u_n)  | < \frac{1}{2}  \left\{  q_{\beta'}  \Phi^{[k-1]} +  \Phi^{[k]}
\right\} (u_n) \\
& \hspace{.5cm} \text{ and } 2\leq k \leq N-1, \\
s_1 &\text{if } |\hat \Phi^p_n (u_n) | \geq \frac{1}{2} \left\{ q_{\beta'} \Phi^{[1]}
  +\Phi^{[2]} \right\} (u_n),\\
s_N  &\text{if }  |\hat  \Phi^p_n  (u_n) |  <  \frac{1}{2} \left\{  q_{\beta'}
  \Phi^{[N-1]} +\Phi^{[N]} \right\} (u_n),\\
\end{array}
\right.
\end{equation}
where  $\{u_n\}_{n \geq  0}$ is  a  sequence of  positive  real numbers
growing  to infinity  and  to be  chosen  later. See  Figure~\ref{fig} for  an
illustration of this procedure.

This  estimation procedure is  well-defined for  large enough  $n$ as  for any
$2\leq k \leq N-1$, we have $\{ q_{\beta'} \Phi^{[k]} +\Phi^{[k+1]} \} (u_n) \leq
\{  q_{\beta'}  \Phi^{[k-1]} +  \Phi^{[k]}\} (u_n) $

This procedure is proved to be consistent, with an exponential rate of convergence, in the following proposition.

\begin{proposition}\label{conv_s}
Under  assumptions {\bf (S)}  and {\bf (A)}, consider
the estimation procedure given by \eqref{estim_s} where
$$
u_n= \left( \frac{\log n} {2} -\frac{2\beta'+a \smax}{2\smax} \log \log n \right) ^{1/\smax} ,
$$
for some fixed $a >1$
and the equidistant grid $\smin=s_1<s_2<\ldots<s_N=\smax$ is chosen as
$$
|s_{k+1} -s_k| =  d_n =  \smax (\log n)^{-1}(\log \log n)^{-1}
 \quad ; N -1 =(\smax-\smin)/ d_n.
$$
Then, $\hat s_n$ is strongly consistent, i.e.
$$  \lim_{n \to \infty } \hat s_n = s  \; \; ; \; \;  \pr_{f,s} -\text{almost surely}. $$
Moreover, for each  number of observations $n$, denote  by $\sgrid$ the unique
point $s_k$ on the grid such that $s_k \leq s<s_{k+1}$. We have
$$ \pr_{f,s}(\hat s_n \neq \sgrid) \leq \exp\left( -\frac{A^2}{4} (\log n)^a  (1+o(1))\right) ,$$
where $A$ is defined in Assumption {\bf (A)} and $a >1$ depends on the choice of $u_n$.
\end{proposition}

\begin{remark}
The result remains valid for any sequence $d_n$ satisfying
$$
d_n u_n^{\smax} \log u_n \leq 1 \quad \text{and} \quad \log (1/d_n) =o((\log n)^a).
$$
\end{remark}

\subsection{Adaptive estimation and tests}\label{sec:apres_estim_s}
For the rest of this section, we shall assume that the unknown density $f$ belongs to
some Sobolev class $\mathcal{F}(0,0,\beta,L)$
where $\beta >0$ is the smoothness parameter and $L$ is a positive constant.
We assume that the unknown parameter $\beta$ belongs to some known interval $ [\betainf, \betasup]$.

We  now plug  the preliminary  estimator  of $s$  in the  usual estimation  and
testing procedures. 

Let us introduce the kernel deconvolution estimator $\hat K_n$ built on the
preliminary  estimation   of  $s$  and   defined  by  its   Fourier  transform
$\Phi^{\hat K_n}$,
\begin{eqnarray}
  \label{kernel_s_inconnu}
  \Phi^{\hat K_n} (u) &=& \exp \left\{\left(\frac{|u|}{\hat h_n}\right)^{\hat s_n}\right\} 1_{|u| \leq 1}  ,\\
 \text{where }  \hat  h_n &=& \left( \frac{\log n}{2} -\frac{ \betasup- \hat s_n +1/2}{\hat s_n} \log \log n \right)^{-1/ \hat s_n} .
  \label{bandwidth_s_inconnu}
\end{eqnarray}
Note that both the bandwidth sequence $\hat h_n$ and the kernel $\hat K_n$ are
random and depend on  observations $Y_1, \ldots, Y_n$.  Now, the
 estimator of $f$ is given by
\begin{equation}
  \label{deconv_estim}
\hat  f_{n}  (x)  =\frac{1}{n  \hat h_n}  \sum_{j=1}^n  \hat  K_n\left(
  \frac{Y_j -x} {\hat h_n} \right).
\end{equation}
This estimation  procedure is consistent  and adaptively achieves  the minimax
rate of convergence when considering  unknown densities $f$ in the union of Sobolev
balls $\mathcal{F}(0,0,\beta, L)$ with $\beta \in [\betainf, \betasup] \subset (1/2; +\infty)$
and unknown smoothness parameter $s \in [\smin;\smax] $.

Note that  when a function belongs to $\mathcal{F}(0,0,\beta,L)$  and assumption {\bf  (A)}
is  fulfilled, we necessarily have $\beta' > \beta +1/2$.

\begin{corollary}\label{cor:deconv_s_inconnu}
Under  assumptions  {\bf (S)}  and  {\bf  (A)},  for any  $\betasup>\betainf >1/2$,  the
estimation procedure given by  \eqref{deconv_estim} which uses estimator $\hat
s_n$  defined  by  \eqref{estim_s}  with  parameter  values:  $u_n$  given  by
Proposition~\ref{conv_s} with $a >\smax/ \smin$,
$$ d_n = \min \left\{ (\log n)^{-(\betasup -1/2)/ \smin} , \smax(\log n \log \log n)^{-1}\right\},
$$
satisfies, for any real number $x$,
$$
\limsup_{n\to \infty} \sup_{s \in [\smin;\smax] } \sup_{\beta \in [\betainf, \betasup]}
\sup_{f \in \mathcal{F}(0,0,\beta,  L)} (\log n) ^{(2\beta-1)/ s} \esp_{f,s} | \hat f_n(x) -f(x)|^2 < \infty.
$$
Moreover, this rate of convergence is asymptotically adaptive optimal.
\end{corollary}

\begin{remark}
This result  is obtained by using  that, with high  probability, the estimator
$\hat  s_n$ is  equal to  the point  $s_k$ on  the grid  such that  $s_k\leq s
<s_{k+1}$ (see Proposition~\ref{conv_s}). Then, using the deconvolution kernel
built  on $s_k$  is as  good  as using  the true  value  $s$, as  soon as  the
difference $|s_k -s|$  is sufficiently small (which is ensured  by the size of
the grid).  Note that the fact that we underestimate $s$ by using $s_k \leq s$ is
rather important as deconvolution with overestimated $s$ would lead to unbounded risk.
\end{remark}

Note that the optimality of this procedure is a direct consequence
of a result by \citep{Efromovich} where he considers the
convolution model for circular data with $\beta$ and $s$ fixed and
known. Therefore we may say that there is no loss due to adaptation
neither with respect to $s$ or $\beta$.

Using the same kernel estimator \eqref{kernel_s_inconnu} and the same random
bandwidth~\eqref{bandwidth_s_inconnu}, we define
\begin{equation} \label{Testim_s_inconnu}
\hat T_{n}
=\frac{2}{n(n-1)}  \sumsum_{1\leq k  <  j \leq  n}  < \frac  1 {\hat  h_n}\hat
K_{n}\left( \frac  {\cdot -Y_k}{\hat h_n} \right)  \; , \; \frac  1 {\hat h_n}
\hat K_{n}\left( \frac{\cdot -Y_j} {\hat h_n}\right) >.
\end{equation}

\begin{corollary}\label{cor:intf2_s_inconnu}
Under  assumptions  {\bf (S)}  and  {\bf  (A)},  for any  $\betasup>\betainf >0$,  the
estimation  procedure given by  \eqref{Testim_s_inconnu} which  uses estimator
$\hat   s_n$  defined  by   \eqref{estim_s}  with   parameter  values:  $u_n$  given  by
Proposition~\ref{conv_s} with $a >\smax/ \smin$,
$$ d_n = \min \left\{ (\log n)^{-2\betasup/ \smin} , \smax(\log n \log \log n)^{-1}\right\},
$$
satisfies,
$$
\limsup_{n\to \infty} \sup_{s \in [\smin;\smax] } \sup_{\beta \in [\betainf, \betasup]}
\sup_{f  \in \mathcal{F}(0,0,\beta,  L)}  \left( \log  n\right)^{2\beta/s}
\left\{\esp_{f,s} \left|  \hat T_{n  } -\int f^2  \right |^2  \right\}^{1/2} <
\infty.
$$
Moreover, under additional Assumption {\bf (E)}, this rate of convergence is asymptotically adaptive optimal.
\end{corollary}

The rate of convergence of this procedure  is the same as in the case of known
self-similarity  index  $s$ and  known  smoothness  parameter  $\beta$. It  is thus
asymptotically adaptive optimal according to results obtained by \citep{Butucea}.\\

Let us now define, for any $f_0 \in \mathcal{F}(0,0,\betasup, L)$,
\begin{equation} \label{Tstat_s_inconnu}
\hat T_{n}^0
=\frac{2}{n(n-1)} \sumsum_{1\leq k < j \leq n}
< \frac 1 {\hat h_n}\hat K_{n}\left( \frac {\cdot -Y_k}{\hat h_n} \right) -f_0
\; , \; \frac 1 {\hat h_n} \hat K_{n}\left( \frac{\cdot -Y_j} {\hat h_n}\right) -f_0 >.
\end{equation}
This statistic is used for goodness-of-fit testing of the hypothesis
\begin{eqnarray*}
H_0 &: & f = f_0 \nonumber \\
\text{versus   }   H_1(\mathcal{C},   \Psi_n)    &:&   f   \in   \cup_{\beta   \in
  [\betainf,\betasup]} \{  f \in \mathcal{F}(0,0,\beta,L) \text{ and  } \psi_{n,\beta}^{-2}
\|f-f_0\|_2^2 \geq \mathcal{C} \} .
\end{eqnarray*}
The test is constructed as usual
\begin{equation}\label{test_s_inconnu}
\Delta_n^\star =\left\{
\begin{array}{ll}
1 & \text{if }  | \hat T_{n}^0 | \hat t_{n}^{-2} >\mathcal{C}^\star\\
0 & \text{otherwise},
\end{array}
\right.
\end{equation}
for some constant $\mathcal{C}^\star>0$ and a {\bf random} threshold $\hat t_{n}^2 $ to be specified.

\begin{corollary}\label{cor:test_s_inconnu}
Under  assumptions {\bf (S)} and  {\bf (A)}, for any $0 < \betainf< \betasup$, any $L>0$
and for any $f_0 \in \mathcal{F}(0,0,\betasup, L)$,  consider
the  testing procedure  given by  \eqref{test_s_inconnu} which  uses  the test
statistic \eqref{Tstat_s_inconnu} with estimator $\hat s_n$ defined by
\eqref{estim_s} with parameter values:  $u_n$  given  by
Proposition~\ref{conv_s} with $a >1$,
$$ d_n  = \min \left\{  (\log n)^{-\betasup/ \smin}  , \smax(\log n  \log \log
  n)^{-1}\right\} ,
$$
with random threshold and (slightly modified) random bandwidth
$$
\hat t_n^2 = \left( \frac{\log n}{2} \right) ^{-2 \betasup / \hat s_n }
\quad ; \quad
\hat h_n =\left( \frac{\log n}{2} - \frac{2\betasup}{\hat s_n} \log \log n \right)^{-1/\hat s_n}
$$
and any large enough positive constant $\mathcal{C}^\star$.
This testing procedure satisfies \eqref{uptest} for any $\epsilon \in (0,1)$ with testing rate
$$
\Psi_n = \{\psi_{n,\beta}\}_{\beta \in [\betainf, \betasup]} \text{ given by }
\psi_{n, \beta} = \left( \frac{\log n}{2} \right) ^{- \beta /s } .
$$
Moreover, if  $f_0 \in  \mathcal{F}(0,0,\betasup, cL)$ for  some $0<c<1$ and  if Assumptions
{\bf (T)} and {\bf(E)} hold,
then this testing rate is asymptotically adaptive optimal over the family of
classes $\{ \mathcal{F}(0,0,\beta, L), \beta \in [\betainf; \betasup] \} $
and for any $s \in [\smin;\smax]$ (i.e. \eqref{lowtest} holds).
\end{corollary}

Adaptive optimality  (namely  \eqref{lowtest}) of  this testing  procedure
directly follows from \citep{Butucea} as there  is no loss due to adaptation to
$\beta$ nor to $s$.  Note also that the case of known  $s$ and adaptation only
with respect to $\beta$ is included in our results and is entirely new.

%
\section{Auxiliary result: Berry-Esseen inequality for degenerate $U$-statistics of order 2}
\label{sec:Ustats}
This section is dedicated to  the statement of a non-uniform Berry-Esseen type
theorem for degenerate $U$-statistics.
It  draws its  inspiration from  \citep{Hall}  which provides  a central  limit
theorem for  degenerate $U$-statistics. Given  a sample $Y_1, \ldots,  Y_n$ of
i.i.d. random variables, we shall consider $U$-statistics of the form
$$
U_n =\sumsum_{1\leq i<j\leq n} H(Y_i,Y_j),
$$
where $H$ is a symmetric function.  We may assume, without loss of generality,
that $\esp\{H(Y_1,Y_2)\}  =0$ and  thus $U_n$ is  centered. We shall  focus on
{\bf degenerate} $U$-statistics, namely
$$
\esp\{H(Y_1,Y_2) | Y_1\} =0 \text{ , almost surely.}
$$
Limit theorems for degenerate $U$-statistics when $H$ is fixed (independent of
the sample size $n$)  are well-known and can be found in  any monograph on the
subject (see for instance \citep{BorKor}).
In that  case, the limit distribution  is a linear  combination of independent
and   centered   $\chi^2(1)$  (chi-square   with   one   degree  of   freedom)
distributions. However,  as noticed in \citep{Hall}, a normal distribution may
result in some cases where $H$  depends on $n$. In such a context, \citep{Hall}
provides  a central  limit theorem.  But  this result  is not  enough for  our
purpose (namely, optimality in Theorem~\ref{th:bruit_poly_f0sob}).  Indeed, we need to control the
convergence  to zero  of the  difference between  the  cumulative distribution
function (cdf) of our $U$-statistic, and the cdf of the Gaussian distribution.
Such a result may be derived using classical Martingale methods.

In  the rest  of this  section,  $n$ is  fixed. Denote  by $\mathcal{F}_i$  the
$\sigma$-field generated by the random variables
$\{Y_1, \ldots, Y_i\}$. Define
$$
v_n^2 = \esp(U_n^2) \quad ; \quad Z_i= \frac{1}{v_n} \sum_{j=1}^{i-1} H(Y_i,Y_j) ,
\quad 2 \leq i \leq n
$$
and note that  as the $U$-statistic is degenerate,  we have $\esp(Z_i |Y_1,
\ldots, Y_{i-1}) =0$. Thus,
$$
S_k = \sum_{i=2}^k Z_i , \quad 2\leq k \leq n,
$$
is a centered Martingale (with respect to the filtration $\{\mathcal{F}_k\}_{k\geq 2}$)
and $S_n = v_n^{-1} U_n$.
We  use    a non-uniform  Berry-Esseen    type    theorem   for
Martingales provided by \citep{HallHeyde}, Theorem 3.9.
Denote by $\phi$ the cdf of the standard Normal distribution and introduce the
conditional variance of the increments $Z_j$'s,
$$
V_n^2    =\sum_{i=2}^n    \esp(Z_i^2   |\mathcal{F}_{i-1})    =\frac{1}{v_n^2}
\sum_{i=2}^n  \esp  \left\{  \left(  \sum_{j=1}^{i-1} H(Y_i,Y_j)  \right)^2  \Big|
  \mathcal{F}_{i-1} \right\}.
$$

\begin{theorem}\label{th:Ustats}
Fix $0 < \delta \leq 1$ and define
$$
L_n = \sum_{i=2}^n \esp |Z_i|^{2+2\delta} + \esp|V_n^2-1|^{1+\delta} .
$$
There exists a positive constant $C$  (depending only on $\delta$)
such that for any $0<\epsilon<1/2$ and
any real $x$
\begin{equation*}
|\pr   (   U_n   \leq   x)   -\phi(x/v_n)   |   \leq   16   \epsilon^{1/2}
\exp\left(-\frac{x^2}{4v_n^2}\right)+ \frac{C}{\epsilon^{1+\delta }} L_n .
\end{equation*}
\end{theorem}

%
\section{Proofs}\label{sec:proofs}
We use $C$ to denote an absolute constant
which values may change along the lines.

\begin{proof} [Proof of Theorem~\ref{th:bruit_poly_f0sob} (Upper bound)]
Let us give the sketch of proof concerning the upper-bound of the test.
The  statistic $T_{n, h^i}$  will  be abbreviated  by  $T_{n,i}$.
We first need to control the first-type error of the test.
\begin{eqnarray*}
\pr_0(\Delta_n^\star =1) &=&\pr_0(\exists i\in \{1,  \ldots , N+1 \} \text{ such
  that } |T_{n,i}| > \mathcal{C}^\star t_{n,i}^2 )\\
&\leq  &  \sum_{i=1}^{N+1}  \pr_0(|T_{n,i} -\esp_0(T_{n,i})|  >  \mathcal{C}^\star
t_{n,i}^2 -\esp_0(T_{n,i})).
\end{eqnarray*}
The proof relies on the two following lemmas.

\begin{lemma}\label{lem1}
For any large enough  $\mathcal{C}^\star >0$, we have
$$
\sum_{i=1}^N \pr_0(|T_{n,i} -\esp_0(T_{n,i})|  >  \mathcal{C}^\star t_{n,i}^2 -\esp_0(T_{n,i})) =o(1).
$$
\end{lemma}

\begin{lemma}\label{lem2}
For  large enough  $\mathcal{C}^\star$, there is some
$\epsilon \in  (0,1)$, such that
$$
 \pr_0(|T_{n,N+1} -\esp_0(T_{n,N+1})|  >  \mathcal{C}^\star
t_{n,N+1}^2 -\esp_0(T_{n,N+1})) \leq \epsilon.
$$
\end{lemma}

Lemma~\ref{lem1}     relies    on     the     Berry-Esseen    type     theorem
(Theorem~\ref{th:Ustats})
presented in Section~\ref{sec:Ustats}. Its proof is postponed
to the  very end of the present  proof. Proof of Lemma~\ref{lem2}  is easy and
omitted. {\it Note for the referee: omitted proofs appear in the appendix.}

Thus, the first type error term is as small as we need, as soon as we choose a large
enough constant  $\mathcal{C}^\star >0$ in  \eqref{test}.
We  now  focus  on the  second-type  error  of  the  test. We write
\begin{multline*}
\sup_{\tau     \in    \mathcal{T}}     \sup_{f     \in    \mathcal{F}(\tau,L)}
\pr_f(\Delta_n^\star =0) \\
\leq 1_{\rinf >0} \sup_{ r \in [\rinf;\rsup], \alpha\geq \alphamin, \beta\in [\betainf,\betasup]}
\sup_{\substack{f \in  \mathcal{F}(\tau,L)\\ \|f-f_0\|_2^2  \geq
    \mathcal{C}\psi^2_{n,\tau}}}
\pr_f ( | T_{n,N+1}| \leq \mathcal{C}^\star t_{n,N+1}^2) \\
+ 1_{\rinf =\rsup=0} \sup_{\alpha\geq \alphamin, \beta \in [\betainf;\betasup]}
\sup_{\substack{f \in  \mathcal{F}(\alpha,0,\beta,L)\\ \|f-f_0\|_2^2  \geq
    \mathcal{C}\psi^2_{n,(\alpha,0,\beta)}}}
\pr_f ( \forall 1\leq i \leq N, \; | T_{n,i}| \leq \mathcal{C}^\star t_{n,i}^2) .
\end{multline*}
Note that when the function $f$ in the alternative is supersmooth
($\rinf>0$), we only need the last test (with index $N+1$), whereas
when it is ordinary smooth ($\rinf=\rsup=0$), we use the family of
tests with indexes $i\leq N$. In this second case, we use in fact
only the test based on parameter $\beta_f$ defined as  the smallest
point on  the grid  larger than  $\beta$ (see  the  proof of
Lemma~\ref{lem4} below).
\begin{lemma}\label{lem4}
We have
$$
\sup_{\alpha \geq   \alphamin}      \sup_{\beta     \in     [\betainf;\betasup]}
\sup_{\substack{f  \in \mathcal{F}(\alpha,  0,  \beta,L)\\ \|f-f_0\|_2^2  \geq
    \mathcal{C}\psi^2_{n,(\alpha,0,\beta)}}}
\pr_f ( \forall 1\leq i \leq N, \; | T_{n,i}| \leq \mathcal{C}^\star t_{n,i}^2) =o(1).
$$
\end{lemma}

\begin{lemma}\label{lem3}
Fix $\rinf>0$,  for any $\alpha  \geq \alphamin, r\in [\rinf;\rsup],  \beta \in
[\betainf;\betasup]$. For  any $\epsilon \in  (0;1)$, there exists  some large
enough $\mathcal{C}^0$ such that for any $\mathcal{C}>\mathcal{C}^0$ and any $f \in \mathcal{F}(\alpha, r, \beta,L)$ such that
$\|f-f_0\|_2^2 \geq \mathcal{C}\psi_{n,(\alpha,r,\beta)}$, we have
$$
 \pr_f ( | T_{n,N+1}| \leq \mathcal{C}^\star t_{n,N+1}^2) \leq \epsilon.
$$
\end{lemma}

The proof of Lemma~\ref{lem4} (resp. \ref{lem3}) is postponed (resp.
omitted) to the very end of the present proof. Thus, the second type
error of the test converges to zero. This ends the proof of
\eqref{uptest}.
\end{proof}
\\

We now present the proofs of the lemmas.

\begin{proof}[Proof of Lemma~\ref{lem1}]
Let us set $\rho_n =(\log
\log  n )^{-1/2}$  and  fix $1\leq  i \leq  N$.  We use  the obvious  notation
$p_0=f_0 * g$.
As we have
\begin{eqnarray*}
\esp_0(T_{n,i})= \| K_{h^i}*p_0 -f_0\|_2^2= \|J_{h^i} *f_0 -f_0\|_2^2, \\
\text{and } <K_h(\cdot -Y_1) -J_h*f_0 , J_h*f_0 -f_0 > =0
\end{eqnarray*}
we easily get
$$
T_{n,i} -\esp_0(T_{n,i}) =\frac{2}{n(n-1)}  \sumsum_{1\leq k<  j  \leq n}
<K_{h^i}(\cdot -Y_k) -J_{h^i}*f_0,K_{h^i}(\cdot -Y_j) -J_{h^i}*f_0>.
$$
Let us set
$$H(Y_j,Y_k)   =    2\{n(n-1)\}^{-1}   <K_{h^i}(\cdot   -Y_k)
-J_{h^i}*f_0,K_{h^i}(\cdot -Y_j) -J_{h^i}*f_0>
$$
and note that  $H$ is a symmetric function  with $\esp_0 \{H(Y_1,Y_2)\}=0$ and
\\
$\esp_0 \{H(Y_1, Y_2) |Y_1\} =0$.  As a consequence, $T_{n,i}
-\esp_0(T_{n,i} )$ is a degenerate $U$-statistic.  Using
Theorem~\ref{th:Ustats}  (and the  notation  of Section~\ref{sec:Ustats})  to
control its cdf, we  get that for any $0< \delta \leq  1$, for any $0<\varepsilon
<1/2$ and any $x$
\begin{multline*}
|\pr_0 ( T_{n,i}-\esp_0(T_{n,i}) > x) -(1-\phi(x/v_n))  | \\
\leq   16   \varepsilon^{1/2}   \exp\left(-\frac{x^2}{4v_n^2}\right)+
\frac{C}{\varepsilon^{1+\delta}} \left  \{ \sum_{i=2}^n \esp_0  |Z_i|^{2+2\delta} +
  \esp_0|V_n^2-1|^{1+\delta} \right \} ,
\end{multline*}
where $v_n^2 = \textrm{Var}_0(T_{n,i})$ and
$$
Z_i = \frac{1}{v_n} \sum_{j=1}^{i-1} H(Y_i,Y_j) \quad \mbox{ and }
V_n^2 = \sum_{i=2}^n \esp_0 (Z_i^2  | \mathcal{F}_{i-1})
$$
as in Section~\ref{sec:Ustats}. Choose  $\delta =1$  and consider
$\varepsilon$ as  a constant  (optimization in $\varepsilon$ is not
necessary in our context), thus
\begin{multline} \label{control_x}
|\pr_0 ( T_{n,i}-\esp_0(T_{n,i}) > x) -(1-\phi(x/v_n))  | \\
\leq C \exp\left(-\frac{x^2}{4v_n^2}\right)+ C \left \{ \sum_{i=2}^n \esp_0 |Z_i|^{4} +
\esp_0|V_n^2-1|^{2} \right \} .
\end{multline}
We want to apply this inequality at the point $x= \mathcal{C}^\star t_{n,i}^2 -\esp_0(T_{n,i})$.
First, note that
$$
\esp_0(T_{n,i}) =\|J_{h^i} *  f_0 -f_0\|_2^2 =\frac{1}{2\pi} \int_{|u| >
  1/(h^i)} |\Phi_0(u)|^{2} du \leq L (h^i)^{2\betasup} \leq L t_{n,i}^2 ,
$$
leading to
$$
x \geq (\mathcal{C}^\star -L)  t_{n,i}^2
= (\mathcal{C}^\star -L)(n\rho_n)^{-4\beta_i /(4\beta_i +4\sigma +1)}
$$
and we choose $\mathcal{C}^\star >L$.
Now, the variance term $v_n^2$ satisfies (see \citep{Butucea})
$$
v_n^2 = \esp_0 (T_{n,i} -\esp_0(T_{n,i}))^2 =\frac{C}{n^2 (h^i)^{4\sigma+1}}
(1+o(1)).
$$
Using the choice of  the bandwidth $h^i$, we obtain a bound  of the first term
in \eqref{control_x}
$$
C \exp\left(-\frac{x^2}{4v_n^2}\right) \leq C
  \exp \left( -\frac{(\mathcal{C}^\star)^2} {C'} \rho_n ^{-2}\right) = C (\log
  n )^{- b},
$$
where $b = (\mathcal{C}^\star) ^2 /(C')$ can be chosen as large as we need.
Let us  deal with  the other terms  appearing in \eqref{control_x}.  For large
enough $n$,
\begin{multline*}
|<K_{h^i}(\cdot -Y_k) -J_{h^i}*f_0,K_{h^i}(\cdot -Y_j) -J_{h^i}*f_0>|\\
 \leq  \frac  {2}{\pi}  \int_{|u|  \leq  1/h^i}  |\Phi^g(u)|^{-2}  du  \leq  \frac{C}{
 (h^i)^{2\sigma +1}}
 \end{multline*}
and thus, for any $p\geq 2$,
\begin{equation*}
\esp_0\{ |H(Y_1, Y_2)|^{2p}\} \leq C n^{-4p} (h^i)^{-2p(2\sigma +1)} .
\end{equation*}
This leads to
 \begin{multline*}
 \sum_{i=2}^n \esp_0 |Z_i|^{4}  \leq  \frac{1}{v_n^4} \sum_{i=2}^{n}
 \left(\sum_{j=1}^{i-1} \esp_0(H(Y_i,Y_j)^4) +  3 \sumsum_{1 \leq j  \neq k \leq
     i-1} \esp_0(H(Y_i,Y_j)^2 H(Y_i,Y_k)^2)
 \right)\\
   \leq  \frac{1}{v_n^4} \sum_{i=2}^{n}
 \left((i-1) \esp_0 (H(Y_1,Y_2)^4) + 3 (i-1)(i-2) \esp_0 (H(Y_1,Y_2)^2 H(Y_1,Y_3)^2)
 \right)\\
  \leq  \frac{O(1)}{v_n^4} n^2 \esp_0(H(Y_1,Y_2)^4) + \frac{O(1)}{v_n^4}
 n^3 \esp_0(H(Y_1,Y_2)^2 H(Y_1,Y_3)^2)\\
  \leq  O(1)\frac{n^3}{n^8 (h^i)^{4(2\sigma + 1)}} n^4 (h^i) ^{2(4\sigma +1)}
 = \frac{O(1)}{n(h^i)^2}.
\end{multline*}

Moreover, following the lines of the proof of Theorem~1 in  \citep{Hall} we get
$$
\esp_0 |V_n^2 -1| ^{2}  \leq \frac{1}{v_n^4} \left(\esp_0(G^2(Y_1,Y_2))+
\frac{1}{n}\esp_0(H^4 (Y_1,Y_2)) \right)  ,
$$
where $G(x,y) =\esp_0(H(Y_1,x) H(Y_1,y)) $. In \citep{ButuceaSORT} this last term was
bounded from above for this model by $C h^i$ so
$$
\esp_0 |V_n^2 -1| ^{2}  \leq C h^i.
$$

Returning  to  \eqref{control_x} we  finally  get  for $x=
\mathcal{C}^\star  t_{n,i}^2 -\esp_0(T_{n,i})$,
\begin{eqnarray*}
|\pr_0(T_{n,i} -\esp_0(T_{n,i}) >x) -\{1-\phi(x/v_n)\} |
 \leq  C\left( (\log n )^{-b} +  h^i \right) \leq  C(\log n)^{-b} .
 \end{eqnarray*}
 Finally we obtain, for $b$ large when $\mathcal{C}^\star$ is
 large
\begin{multline*}
\sum_{i=1}^N \pr_0(|T_{n,i} -\esp_0(T_{n,i})|  >  \mathcal{C}^\star t_{n,i}^2 -\esp_0(T_{n,i}))
     \leq  N ( 1-\phi(x/v_n ) +C(\log n)^{-b}) \\
\leq  C N \left( v_n x^{-1} \exp(-x^2/(2v_n^2)) + (\log
n)^{-b}\right) \leq  C N \rho_n (\log n)^{-b} \leq C \frac{(\log
\log n)^{-1/2} }{\log n^{b-1}}.
\end{multline*}
\end{proof}

\begin{proof} [Proof of Lemma~\ref{lem4}]
 When $\rsup=\rinf=0$, let us fix  some  constant
$\mathcal{C}>\mathcal{C}^0  $ ($\mathcal{C}^0  $ will be chosen later)
and  a   density  $f$  belonging  to  $\mathcal{F}(\alpha,0,\beta,L)$ for some
unknown $\alpha > \alphamin$ and $\beta  \in [\betainf; \betasup]$ which satisfies
$\|f-f_0\|_2^2 \geq \mathcal{C} \psi^2_{n,(\alpha,0,\beta)}$ (choose $\beta$
as the largest one). In this proof, we abbreviate $\psi_{n,(\alpha,0,\beta)}$ to
$\psi_{n,\beta}$ since in this case, the rate only depends on $\beta$.
We define $\beta_f$ as
the smallest point  on the finite grid $\{\betainf =  \beta_1 < \beta_2 <\ldots
<\beta_N =\betasup\}$ such that $\beta \leq \beta_f$
\begin{multline}\label{betaf}
\beta_f \in \{\betainf  = \beta_0 < \beta_1 <\ldots  <\beta_N =\betasup\} , \;
\; f \in \mathcal{F}(\alpha, 0,\beta, L), \|f-f_0\|_2^2 \geq \mathcal{C}\psi^2_{n,\beta} ,\\
\beta \leq  \beta_f \text{ and } \forall  \beta_i < \beta_f, \text{  we have }
\beta >\beta_i .
\end{multline}
We shall abbreviate  to $h_f$, $t_{n,f}^2$ and $T_{n,f}$  the bandwidth, the
threshold  (both  defined   in  Theorem~\ref{th:bruit_poly_f0sob})  and  the  statistic
\eqref{Tstat} corresponding to parameter $\beta_f$. We write
\begin{eqnarray}
 &&\pr_f (\forall i \in  \{1, \ldots , N\}  , \; \;
| T_{n,i}| \leq \mathcal{C}^\star t_{n,i}^2) \nonumber \\
&\leq &  \pr_f (|T_{n,f} -\esp_f(T_{n,f})|  \geq -\mathcal{C}^\star t_{n,f}^2
+ \esp_f(T_{n,f}) )  \nonumber \\
&\leq & \pr_f (|T_{n,f} -\esp_f(T_{n,f})| \geq \|f-f_0\|_2^2
- \mathcal{C}^\star t_{n,f}^2  + B_f(T_{n,f})) ,\label{prem0}
\end{eqnarray}
where
$$
B_f(T_{n,f}) = \esp_f(T_{n,f}) - \|f-f_0\|_2^2= \|J_h *f\|_2^2 -\|f\|_2^2
+2 \langle f- J_h * f, f_0   \rangle
$$
is in fact a bias term. It satisfies
\begin{eqnarray*}
|B_f(T_{n,f})| &\leq & \int_{|u|\geq 1/h_f} |\Phi (u)|^2 du
+2 (\int_{|u|\geq 1/h_f} |\Phi (u)|^2 du \int_{|u|\geq 1/h_f} |\Phi_0 (u)|^2 du)^{1/2}\\
&\leq & L e^{-2\alphamin}(h_f ^{2\beta}+2 h_f^{\betasup+\beta})
\leq 3e^{-2\alphamin}L h_f^{2\beta},
\end{eqnarray*}
as $f$ belongs to $\mathcal{F}(\alpha, 0 , \beta, L)\subseteq \mathcal{F}(\alphamin, 0 , \beta, L)$.

Let us study the variance term $ \esp_f (T_{n,f}
-\esp_f(T_{n,f}))^2$.
According to \citep{Butucea}, this term is upper-bounded by $w_{n,f}^2$ given by
$$
 \esp_f (T_{n,f} -\esp_f(T_{n,f}))^2
\leq \frac{C }{n^2 h_f^{4\sigma +1}} +\frac{4 \Omega_g^2(f-f_0)} {n} 1_{\beta\geq \sigma}
= w_{n,f}^2,
$$
and $\Omega_g(f-f_0)$  is a constant  depending on $f$  and $g$ (but not  $n$) and
satisfying $|\Omega_g^2(f-f_0) | \leq C\|f- f_0\|_2^{2-2\sigma/\beta}$ (see proof of
Theorem 6 in \citep{Butucea}).

Using  Markov's  inequality,  this  leads  to the  following  upper  bound  of
(\ref{prem0})
\begin{eqnarray*}
   \frac{w_{n,f}^2}
   {( \|f-f_0\|_2 ^2 -\mathcal{C}^\star t_{n,f}^2- 3e^{-2\alphamin}L h_f^{2\beta})^2} .
\end{eqnarray*}

We will  proceed differently when  $\beta<\sigma$ and when  $\beta \geq \sigma$.  Let us
first consider the term concerning $\beta <\sigma$.
The point is to use that $f $ satisfies $\|f-f_0\|_2^2 \geq
\mathcal{C}\psi_{n,\beta}^2$.
Note that we have $\beta_f \geq \beta$, constants $\mathcal{C} >
\mathcal{C}^\star$ and
\begin{eqnarray*}
  \psi_{n,\beta}    ^2   t_{n,f}^{-2}    &=&    (n\rho_n)^{4(\beta_f   -\beta)
 (4\sigma+1)/\{(4\beta_f +4\sigma+1)(4\beta +4\sigma +1)\}},
\end{eqnarray*}
ensuring that the term $\mathcal{C}\psi_{n,\beta}^2 - \mathcal{C}^\star
t_{n,f}^2 $ is always positive.  Moreover, as $0 \geq \beta -\beta_f \geq
- (\betasup-\betainf) / \log n$, we have
\begin{eqnarray*}
  \psi_{n,\beta}^2 h_f^{-2\beta} &=& \exp\left \{ \frac{ 16\beta (\beta -\beta_f) }{
    (4\beta_f +4\sigma+1)(4\beta +4\sigma +1)} \log(n\rho_n) \right \}\\
&\geq  &  \exp   \left\{-  \frac{16\betasup  (\betasup-\betainf)  }{(4\betainf
    +4\sigma+1)^2 } (1+o(1)) \right\} =: \mathcal{C}_1 (1+o(1)) .
\end{eqnarray*}
Thus, we choose $\mathcal{C}^0=\mathcal{C}^\star+3e^{-2\alphamin}L/\mathcal{C}_1$ such that
for any $\mathcal{C}>\mathcal{C}^0$, we have
\begin{eqnarray*}
\|f-f_0\|_2 ^2 -\mathcal{C}^\star t_{n,f}^2- 3e^{-2\alphamin} L h_f^{2\beta}
  \geq  (\mathcal{C}   -
  \mathcal{C}^* - 3e^{-2\alphamin}L/ \mathcal{C}_1 )\psi_{n,\beta}^2 = a\psi_{n,\beta}^2 ,
\end{eqnarray*}
with $a>0$. Thus, we get
\textsc{}\begin{multline*}
\sup_{\alpha>     \alphamin}      \sup_{\beta     \in     [\betainf;\betasup]}
\sup_{\substack{f \in \mathcal{F}(\alpha, 0, \beta,L)\\ \|f-f_0\|_2^2 \geq
 \mathcal{C} \psi^2_{n,\beta}}}
\pr_f (\forall i \in  \{1, \ldots , N\}  , \; \; | T_{n,i}| \leq \mathcal{C}^\star t_{n,i}^2) \\
\leq
\max \left\{ \sup_{\beta<\sigma} \sup_{f }
\frac{C}{n^2    h_f^{4\sigma+1}\psi_{n,\beta}      ^4} ;
\sup_{\beta \geq \sigma}\sup_{f }
  \frac{C\|f-f_0\|_2^{2-2\sigma/\beta}}
  {n ( \|f-f_0\|_2 ^2 -\mathcal{C}^\star t_{n,f}^2-3e^{-2\alphamin} L h_f^{2\beta})^2
   }\right\}.
\end{multline*}

Finally, this leads to the bound
\begin{eqnarray*}
 &&\max \left\{ \sup_{\beta < \sigma} \sup_f \frac{C}{n^2
h_f^{4\sigma+1}\psi_{n,\beta}      ^4}      ;      \sup_{\beta     \geq      \sigma}
  \frac{C } {n \|f-f_0\|_2 ^{2+2\sigma/\beta }(a/\mathcal{C})^2}\right\} \\
  &\leq&
 \max \left\{ \sup_{\beta < \sigma} \sup_f \frac{C}{n^2
h_f^{4\sigma+1}\psi_{n,\beta}      ^4}      ;      \sup_{\beta     \geq      \sigma}
  \frac{C } {n \psi_{n,\beta} ^{2+2\sigma/\beta }}\right\} \leq \rho_n ,
\end{eqnarray*}
which converges to zero as $n$ tends to infinity.
\end{proof}
\\

\begin{proof} [Proof of Theorem~\ref{th:bruit_poly_f0sob} (Lower bound)]

As we already noted after the theorem statement, our test procedure
attains the minimax rate associated to the class
$\mathcal{F}(\alpha_0,0,\betasup,  L)$  where   $f_0$  belongs,  whenever  the
alternative $f$ belongs to classes of
functions smoother than $f_0$. Therefore, the lower bound we need to
prove concerns the optimality of the loss of order $(\log \log
n)^{1/2}$ due to alternatives less smooth than $f_0$.

More precisely, we prove \eqref{lowtest}, where the alternative
$H_1(\mathcal{C}, \Psi_{n})$ is now restricted to
$\cup_{\beta \in [\betainf,\betasup]} \{f  \in \mathcal{F}(0,0,\beta,L) \text{
and } \psi_{n,\beta}^{-2} \|f-f_0\|_2^2 \geq \mathcal{C} \}$ and
$\psi_{n,\beta}$ denotes the rate $\psi_{n,\tau}$ when
$\tau=(0,0,\beta,L)$.

The  general  approach for  proving  such a lower  bound
\eqref{lowtest}  is  to exhibit  a finite  number of  regularities
$\{\beta_k\}_{1\leq k  \leq K}$  and corresponding probability
distributions $\{\pi_k\}_{1\leq  k \leq K}$  on the alternatives
$H_1(\mathcal{C},\psi_{n,\beta_k})$  (more exactly, on parametric
subsets of these alternatives) such that the distance between the
distributions induced by $f_0$ (the density being tested) and the
mean distribution of the alternatives  is small.

\noindent We use a finite grid
$  \bar{\mathcal{B}}  =\{  \beta_1  <\beta_2  <\ldots <  \beta_{K}  \}  \subset
[\betainf,  \betasup]$  such that
$$ \forall  \beta  \in [\betainf,  \betasup],
\exists k : |\beta_{k}-\beta| \leq \frac 1 {\log n}.$$
 To each point $\beta$
in this grid, we associate a bandwidth
$$
h_\beta = \left(  n \rho_n \right)^{-  \frac 2 {4\beta+4\sigma +1}}  , \rho_n = (\log
\log n)^{-1/2}, \quad \text{and } \quad M_\beta = h_\beta ^{-1}.
$$

\noindent We  use the  same deconvolution  kernel as  in \citep{Butucea},  constructed as
follows. Let  $G$ be  defined as  in Lemma 2  in \citep{Butucea}.
The function $G$ is an infinitely differentiable function, compactly supported
on $\displaystyle{\left[-1, \; 0 \right]}$ and such that $\int G=0$.
Then, the deconvolution kernel $H_\beta$  is defined via its Fourier transform
$\Phi^{H_\beta}$ by
$$
\Phi^{H_\beta}(u ) = \Phi^G(h_\beta u) (\Phi^g(u))^{-1} .
$$
Note that the factor $\rho_n$ in the bandwidth's expression corresponds to the
loss for adaptation.

\noindent  We  also  consider   for  each  $\beta$,  a   probability
distribution $\pi_\beta$  (also denoted $\pi_k$  when $\beta=\beta_k$) defined
on $\{-1,+1\}^{M_{\beta}}$ which is in fact the product of Rademacher
distributions on $\{-1,+1\}$ and a parametric subset of $H_1(\mathcal{C},
\psi_{n,\beta})$ containing the following functions
$$
f_{\theta, \beta} (x) = f_0(x)
+ \sum_{j=1}^{M_\beta} \theta_j h_\beta^{\beta+\sigma +1} H_\beta \left( x-x_{j,\beta}   \right) ,
\quad
\left\{
\begin{array}{l}
\theta_j \text{ i.i.d. with } \pr(\theta_j =\pm 1)= 1/2, \\
x_{j,\beta} =j h_\beta \in [0,1].
\end{array}
  \right.
$$
Convolution  of these  functions with  $g$ induces  another parametric  set of
functions
$$ p_{\theta, \beta} (y) =
p_0 (y) + \sum_{j=1}^{M_\beta} \theta_j h_\beta^{\beta+\sigma +1}
 G_\beta \left( y-x_{j,\beta}   \right)$$
where $\displaystyle{ G_\beta (y) = h_\beta^{-1} G\left( y /
 {h_\beta} \right) = H_\beta *  g (y)}$.

\noindent As  established in \citep{Butucea} (Lemmas  2 and 4), for any  $\beta $, any
$\theta \in \{-1,+1\}^{M_\beta}$ and small enough $h_\beta$ (i.e. large enough
$n$) the function $f_{\theta, \beta}$ is a probability density
and belongs to the Sobolev class $\mathcal{F}(0,0,\beta,L)$ and $p_{\theta,\beta}$ is also a
probability density.
Moreover we have
$$ \frac 1 K \sum_{\beta \in \bar{\mathcal B}}
\pi_\beta \left( \| f_{\theta,\beta} - f_0 \|_2^2 \geq  \mathcal{C} \psi_{n,\beta}^2 \right)
\mathop{\longrightarrow}_{n \rightarrow +\infty} 1,$$
which means that for each $\beta $, the random parametric
family  $\{f_{\theta,\beta}\}_{\theta}$ belongs almost surely  (with  respect to  the
measure    $\pi_\beta$)    to    the   alternative    set    $H_1(\mathcal{C},
\psi_{n,\beta})$.
The subset  of functions  which are not  in the  alternative $H_1(\mathcal{C},
\Psi_{n})$ is asymptotically negligible.
 We then have,
\begin{eqnarray*}
 \gamma_n &\triangleq& \inf_{\Delta_n} \left \{ \pr_0 (\Delta_n =1) + \sup_{f \in
     H_1(\mathcal{C},\Psi_n)} \pr_f (\Delta_n =0) \right \}  \\
&\geq  &   \inf_{\Delta_n}  \left  \{   \pr_0  (\Delta_n  =1)   +  \frac{1}{K}
  \sum_{k=1}^K \sup_{f \in  H_1(\mathcal{C},\psi_{n,\beta_k})} \pr_f (\Delta_n =0) \right
\}  \\
&\geq  &   \inf_{\Delta_n}  \left  \{   \pr_0  (\Delta_n  =1)   +  \frac{1}{K}
  \sum_{k=1}^K \Big( \int_{\theta} \pr_{f_{\theta, \beta_k}} (\Delta_n =0) \pi_{k} (d\theta) \right.  \\
&&\hspace{6cm} \left.  -\pi_k(\|f_{\theta,\beta_k} -f_0\|_2^2 <\mathcal{C}\psi_{n,\beta_k}^2) \Big) \right
\} \\
&\geq  &   \inf_{\Delta_n}  \left  \{   \pr_0  (\Delta_n  =1)   +  \frac{1}{K}
  \sum_{k=1}^K \left( \int_{\theta} \pr_{f_{\theta, \beta_k}} (\Delta_n =0) \pi_{k} (d\theta)  \right) \right\}
+o(1) .
\end{eqnarray*}
Let us denote by
\begin{equation*}
  \pi=   \frac{1}{K}   \sum_{k=1}^K   \pi_{k}  \quad \text{
    and } \quad \pr_\pi =\frac{1}{K} \sum_{k=1}^K \pr_k
    =\frac{1}{K} \sum_{k=1}^K
\int_{\theta} \pr_{f_{\theta, \beta_k}} \; \pi_{k} (d \theta) .
\end{equation*}
Those notations lead to
\begin{eqnarray}
\label{minogamma}
\gamma_n & \geq &
\inf_{\Delta_n} \{ \pr_0(\Delta_n =1) + \pr_{\pi} (\Delta_n =0)\}\nonumber \\
& \geq & \inf_{\Delta_n}\left \{ 1 -\int_{\Delta_n =0} d\mathbb{P}_0 +
\int_{\Delta_n =0} d \mathbb{P}_\pi \right\}
\geq 1 - \sup_A \int_A (d\mathbb{P}_0 - d\mathbb{P}_\pi )\nonumber \\
& \geq & 1-\frac{1}{2} \|\pr_\pi -\pr_0\|_1 ,
\end{eqnarray}
where we used Scheff\'e's Lemma.

\noindent The finite grid $\bar{\mathcal{B}}$ is split into subsets
$\bar{\mathcal{B}} =\cup_{l} \bar { \mathcal{B}}_l$ with $\bar{\mathcal{B}}_l \cap
\bar {\mathcal{B}}_k = \emptyset$ when $l \neq k$ and such that
$$
\forall l , \; \; \forall \beta_1 \neq \beta_2 \in \bar {\mathcal{B}}_l, \; \;
\frac{c \log \log n}{\log n} \leq |\beta_1 -\beta_2| .
$$
The number of subsets $\bar {\mathcal{B}}_l$ is denoted by $K_1 =O(\log \log n)$
and   the  cardinality   $|\bar  {\mathcal{B}}_l|$   of  each   subset  $\bar{
  \mathcal{B}}_l$ is  of the  order $O(\log n  /\log \log n)$,  uniformly with
respect to $l$.

\noindent The lower bound \eqref{lowtest} is then obtained from \eqref{minogamma} in the following way
\begin{equation*}
 \gamma_n
\geq 1-\frac{1}{2K_1} \sum_{l=1}^{K_1} \left \| \frac{1}{|\bar {\mathcal{B}}_l|}\sum_{\beta \in \bar
{\mathcal{B}}_l } \pr_{\beta} -\pr_0 \right \|_1 ,
\end{equation*}
where $\pr_\beta = \int_{\theta}  \pr_{f_{\theta, \beta}} \pi_{\beta}(d \theta)$ .

\noindent  Here we do  not want to apply  the triangular inequality to the  whole set of
 indexes  $\bar{\mathcal{B}}$. Indeed, this  would lead  to a  lower bound  equal to
 0. Yet, if we do not apply some sort of triangular inequality, we cannot
deal with  the sum because of too  much dependency. This is  why we introduced
the subsets $\bar {\mathcal B}_l$ with the property that two points in
the same subset $\bar {\mathcal B}_l$ are far enough away from each other.
This technique was already used in \citep{GayraudPouet} for the discrete
regression model.

Let us denote by $\ell_{\beta}$  the likelihood ratio
$$\ell_\beta = \displaystyle{\frac {d \pr_\beta} {d\pr_0} = \int \frac {d\pr_{f_{\theta,
\beta}}} {d\pr_0} \; \pi_\beta\left( d\theta \right)}.$$
We thus have
\begin{multline*}
 \gamma_n
\geq 1-\frac{1}{2K_1} \sum_{l=1}^{K_1} \int \left( \frac{1}{|\bar{\mathcal{B}}_l |} \sum_{\beta \in \bar
{\mathcal{B}}_l } \ell_{\beta} -1 \right) d\pr_0
=  1-\frac{1}{2K_1} \sum_{l=1}^{K_1} \left\|  \frac{1}{|\bar{\mathcal{B}}_l |}
  \sum_{\beta \in \bar {\mathcal{B}}_l }    \ell_{\beta}-1
\right\|_{\Lset_1(\pr_0)} .
\end{multline*}
\noindent Now   we    use   the    usual   inequality   between    $\mathbb{L}_1$  and
$\mathbb{L}_2$-distances to get that
\begin{equation*}
\gamma_n  \geq  1 -
 \frac 1 {2K_1} \sum_{l=1}^{K_1}
\left\| \frac 1 {|\bar {\mathcal{B}}_l|} \sum_{\beta
\in \bar{\mathcal B}_l} \ell_{\beta} -1 \right\|_{\Lset_2(\pr_0)}
  =   1 -  \frac 1 {2K_1} \sum_{l=1}^{K_1}
\left\{\esp_0
 \left( \frac 1 {|\bar {\mathcal{B}}_l|} \sum_{\beta
\in \bar{\mathcal B}_l} \ell_{\beta} - 1 \right)^2
\right\}^{1/2}.
\end{equation*}

\noindent Let us focus on the expected value appearing in the lower
bound. We have
$$
\esp_0 \left(
\frac 1 {|\bar {\mathcal{B}}_l|} \sum_{\beta
\in \bar{\mathcal B}_l} \ell_{\beta} - 1 \right)^2
= \frac 1 {|\bar {\mathcal{B}}_l|^2} \sum_{\beta \in \bar{\mathcal B}_l}
Q_\beta
+ \frac 1 {|\bar {\mathcal{B}}_l|^2} \mathop{\sum_{\beta, \nu \in \bar{\mathcal B}_l}}_
{\beta \neq \nu} Q_{\beta,\nu},$$
where there are two quantities to evaluate
$$
Q_\beta  = \esp_0  \left(  \left(  \ell_{\beta} -  1  \right)^2 \right)  \quad
\text{and} \quad
Q_{\beta,\nu}  = \esp_0 \left(  \ell_{\beta} \ell_{\nu} - 1 \right).
$$

\noindent The first term $Q_\beta$  is treated as in \citep{Butucea}. It corresponds to the computation of a $\chi^2$-distance between the two models induced by $\pr_\beta$ and $\pr_0$ (see term $\Delta^2$ in \citep{Butucea}). Indeed we
have
$$
Q_\beta \leq C M_\beta n^2 h_\beta^{4\beta+4\sigma +2} \leq C \frac 1 {\rho_n^2}.
$$
This upper bound goes to infinity very slowly. The number of $\beta$'s in each $\displaystyle{
\bar{\mathcal B}_l}$ compensates this behaviour
$$
\frac 1 {|\bar {\mathcal{B}}_l|^2} \sum_{\beta \in \bar{\mathcal B}_l}
Q_\beta \leq  \frac 1 {|\bar  {\mathcal{B}}_l| \rho_n^2} =  O\left(
  \frac {(\log \log n )^2} {\log n} \right) =o(1).
$$

\noindent  The  second  term  is  a  new one  (with  respect  to  non-adaptive
case). As $G$ is compactly supported and  the points $\beta$ and $\nu$
are far  away from each other, we  can prove that this  term is asymptotically
negligible.
Recall the expression of the likelihood ratio for a fixed $\beta$
$$ \ell_\beta = \int \frac {d\pr_{f_{\theta, \beta}}} {d\pr_0} \; \pi_\beta \left( d\theta \right)
= \int \prod_{r=1}^n \left( 1 + \sum_{j=1}^{M_\beta} \theta_{j,\beta}
h_\beta^{\beta+\sigma +1} \frac {G_{\beta}\left( Y_r - x_{j,\beta} \right)}
{p_0\left(Y_r \right)} \right) \pi_\beta \left( d\theta \right).$$
Thus,
\begin{multline*}
 \ell_\beta \ell_\nu
 =
\int \frac {d\pr_{f_{\theta,\beta}}} {d\pr_0} \; \pi_\beta\left( d\theta \right)
\int \frac {d\pr_{f_{\theta,\nu}}} {d\pr_0} \; \pi_\nu\left( d\theta \right) \\
  = \int \prod_{r=1}^n \left( 1 + \sum_{j=1}^{M_\beta} \theta_{j,\beta} \;
h_\beta^{\beta+\sigma +1} \frac {G_{\beta}\left( Y_r -x_{j,\beta} \right)}
{p_0\left( Y_r \right)} \right) \\
\times \left( 1 + \sum_{i=1}^{M_\nu} \theta_{i,\nu}
h_\nu^{\nu+\sigma +1} \frac {G_{\nu}\left( Y_r -x_{i,\nu} \right)}
{p_0\left( Y_r \right)} \right)
  \pi_\beta\left( d\theta_{.,\beta} \right)
 \pi_\nu\left( d\theta_{.,\nu} \right).
\end{multline*}
The random variables $\displaystyle{Y_r}$ are i.i.d. and
$\displaystyle{\esp_0 \left( \frac {G_{\beta}\left( Y_r - x_{j,\beta}
 \right)} {p_0\left(Y_r\right)} \right) = 0}$. Thus we have
\begin{multline*}
\esp_{0}\left( \ell_\beta \ell_\nu \right) =
 \int  \left[1 +
\sum_{j=1}^{M_\beta} \sum_{i=1, i \subset j}^{M_\nu} \theta_{j,\beta}
\theta_{i,\nu} h_\beta^{\beta+\sigma +1} h_\nu^{\nu+\sigma +1} \right.\\\left.
\esp_{0}\left(
  \frac {G_{\beta}\left( Y_1 -x_{j,\beta} \right)G_{\nu}\left( Y_1 -x_{i,\nu}
    \right) } {p_0^2\left( Y_1 \right)}
\right) \right] ^n\pi_\beta\left( d\theta_{.,\beta} \right)
 \pi_\nu\left( d\theta_{.,\nu} \right) .
\end{multline*}
where the second  sum concerns only some indexes $i$,   denoted by
$i  \subset j$. This  notation stands  for the  set of  indexes $i$  such that
$[(i-1)h_\beta;ih_\beta] \cap [(j-1)h_\nu;j h_\nu] \neq \emptyset $.
From now on, we fix $\beta>\nu$.
Denote by $G'$ (resp. $p_0'$) the first derivative of $G$ (resp. $p_0$).
(The density $p_0$ is continuously  differentiable as it is the convolution
product $f_0 * g$ where the noise density $g$ is at least continuously
differentiable).

\begin{lemma}\label{lem5}
For any $\beta> \nu$ and any $(i,j)\in \{1,\ldots, M_\nu\} \times
\{1,\ldots,M_\beta\}$, we have
$$ \esp_{0}\left(
  \frac   {G_{\beta}\left(  Y_1  -x_{j,\beta}   \right)G_{\nu}\left(  Y_1
      -x_{i,\nu} \right)} {p_0 ^2\left( Y_1 \right)}
\right) =  \frac {h_\nu} {h_\beta^2} R_{i,j} ,
$$
where $R_{ij} $ satisfies
$$
|R_{i,j}| \leq (\inf_{[0,1]} p_0)^{-1} \|G\|_\infty \|G'\|_\infty (1+o(1))
$$
and $o(1)$ is uniform with respect to $(i,j)$.
 \end{lemma}

The proof of this lemma is omitted.
\noindent Applying Lemma~\ref{lem5}, we get
$$
Q_{\beta,\nu}  +1 =
\int \left[ 1 +  \sum_{j=1}^{M_\beta}
\mathop{\sum_{i=1, i \subset j}^{M_\nu}}
 \theta_{j,\beta} \theta_{i,\nu}
h_\beta^{\beta+\sigma +1} h_\nu^{\nu+\sigma +1}
 \frac {h_\nu} {(h_\beta)^2} R_{i,j} \right]^n
 \pi_\beta\left( d\theta_{.,\beta} \right)
 \pi_\nu\left( d\theta_{.,\nu} \right)  .
$$

\begin{lemma}\label{lem6}
Let $\displaystyle{U}$ be a real valued random variable such that
$\forall k \in {\mathbb N} , \; \; \esp \left( U^{2k+1} \right) = 0$.
We have, for any integer $n \geq 1$,
$$ \esp \left( 1 + U \right)^{n}
\leq 1 + \sum_{k=1}^{\lfloor \frac n 2 \rfloor} \frac { n^{2k}} {\left(2k\right)!}
 \esp \left( U^{2k} \right),$$
where $\lfloor x \rfloor$ is the largest integer which is smaller than $x$.
\end{lemma}

The proof is obvious and therefore omitted.
Apply Lemma~\ref{lem6} to get the inequality
$$
 Q_{\beta,\nu} \leq \sum_{k=1}^{\lfloor \frac  n 2  \rfloor} \frac
 {n^{2k}} {\left( 2k \right)!} (  h_\beta^{ \beta+\sigma -1} h_\nu^{\nu
 +\sigma +2}
)^{2k}
 \esp_\pi \left( \sum_{j=1}^{M_\beta} \sum_{i=1,i\subset
j}^{M_\nu} \theta_{j,\beta} \theta_{i,\nu} R_{i,j}\right)^{2k}.
$$
But the  $\theta$'s are i.i.d.   Rademacher variables   and the
$R_{i,j}$'s are deterministic, thus
$$
 \esp_\pi \left( \sum_{j=1}^{M_\beta} \sum_{i=1,i\subset
j}^{M_\nu} \theta_{j,\beta} \theta_{i,\nu} R_{i,j} \right)^{2k}= \sum_{1\leq j_1,
\ldots,  j_k  \leq  M_\beta   }  \mathop{\sum_{1\leq  i_1,\ldots  ,  i_k  \leq
  M_\nu}}_{\forall l, i_l \subset j_l} \big(\prod_{l=1}^k R_{i_l,j_l}^{2}\big) .
$$
Using the bound on the $R_{i,j}$ given by Lemma~\ref{lem5},
$$
 \esp_\pi \left( \sum_{j=1}^{M_\beta} \sum_{i=1,i\subset
j}^{M_\nu} \theta_{j,\beta} \theta_{i,\nu} R_{i,j} \right)^{2k} \leq
\Big( ( \inf_{[0,1]}  p_0)^{-1} \|G\|_\infty \|G'\|_\infty (1+o(1)) \Big)^{2k}
h_\nu^k .
$$
Indeed,  each index $j_l$  may take  at most  $M_\beta=h_\beta^{-1}$ different
values but  the constraint $i_l\subset j_l$  implies that each  index $i_l$ is
limited to at most $h_\beta/h_\nu$ different values.
Thus we get
\begin{multline*}
 Q_{\beta,\nu} \leq C
\sum_{k=1}^{\lfloor \frac n 2 \rfloor} \frac {n^{2k}} {\left( 2k \right)!}
 \left( C h_\beta^{\beta+\sigma +1} h_\nu^{\nu+\sigma +1}
 \frac {h_\nu} {h_\beta ^2} \right)^{2k} h_\nu^{-k}\\
  \leq  C \sum_{k=1}^{\lfloor \frac n 2 \rfloor }
 \left(  n^2 h_\beta^{2\beta+2\sigma+1/2} h_\nu^{2\nu+2\sigma +1/2}
 \frac {h_\nu^{ 5/ 2}} {h_\beta ^{ 5/ 2}} \right)^k
  \leq  C\sum_{k=1}^{\lfloor \frac n 2 \rfloor}
 \left(   \frac {h_\nu^{ 5/ 2}} {\rho_n^2 h_\beta ^{ 5/ 2}} \right)^k
  \leq C \frac 1 {\rho_n^2} \frac {h_\nu ^{5/2}}
 {h_\beta ^{ 5/ 2}} .
\end{multline*}
As  $\beta>\nu$ both  belong to  some set  $\bar{ \mathcal{B}}_l$,  we have
$\beta-\nu \geq c (\log \log n)/ (\log  n)$ and according to the choice of the
bandwidths,
$$
  \frac {h_\nu^{5/2}} {h_\beta ^{5/ 2}}
  =  \left( n \rho_n \right)^{- \frac {20 \left( \beta - \nu \right)}
{\left(4 \beta +4\sigma +1\right) \left(4 \nu +4\sigma +1\right)}}
  \leq \exp \Big\{ - \frac{20 \; c \log \log n} {(4\betasup+4\sigma +1)^ 2} (1+o(1))\Big\}
\leq (\log n) ^ {-w},
$$
where the constant $w$ (depending on the constant $\displaystyle{c}$ used
in the construction of the  sets $\bar{\mathcal{B}}_l$) can be tailored to our
need. Therefore
$$  \frac 1 {|\bar{\mathcal{B}}_l|^2} \mathop{\sum_{\beta, \nu \in |\bar{\mathcal B}_l|}}_{\beta \neq \nu}
Q_{\beta, \nu} \leq  \frac C {\rho_n^2 \left( \log n \right)^w}$$
which goes to $\displaystyle{0}$ as n goes to $\displaystyle{
+\infty}$.
We finally obtain the upper bound
$$ \esp_0 \left( \left( \frac 1 {|\bar{\mathcal{B}}_l|} \sum_{\beta
\in |\bar{\mathcal B}_l|} \ell_{\beta} - 1 \right)^2 \right)
\leq O \left( \frac 1 {|\bar{\mathcal{B}}_l| \rho_n^2} \right) +
 O\left( \frac 1 {\rho_n^2 \left( \log n \right)^w} \right) =o(1) ,$$
which leads to
$$
\gamma_n  \geq
 1 - \frac 1 2 \frac 1 {K_1} \sum_{l=1}^{K_1}
\left\{  O\left( \frac 1 {|\bar{\mathcal{B}}_l|\rho_n^2} \right) +
 O\left(  \frac 1  {\rho_n^2  \left( \log  n \right)^c}  \right)\right\}^{1/2}
=1+o(1).$$
\end{proof}
\\



\begin{proof}[Proof of Proposition~\ref{conv_s}]
We fix $\epsilon >0$. Now,
$$\pr_{f,s}(|\hat s_n -s| \geq \epsilon) \leq \pr_{f,s}(\hat s_n \neq \sgrid ) +  \pr_{f,s}(|s -\sgrid| \geq \epsilon).
$$
As $|s-\sgrid| \leq d_n$ which converges to zero, we get that for
large enough $n$, the term $\pr_{f,s}(|s -\sgrid | \geq \epsilon)$
is equal to zero. Let us now consider the term $\pr_{f,s}(\hat s_n
\neq \sgrid ) =\pr_{f,s}(\hat s_n > \sgrid ) + \pr_{f,s}(\hat s_n
<\sgrid )  $. Now, $\sgrid$ is equal to some $s_k$ (using the
labeling among the points of the grid). We have
  \begin{eqnarray*}
   \pr_{f,s}(\hat s_n < s_k)    &=& \sum_{j=1}^{k-1}  \pr_{f,s}(\hat s_n = s_j) \\
&\leq & \sum_{j=1}^{k-1}\pr_{f,s} \left( |\hat \Phi^p_n  (u_n) | \geq  \frac{1}{2} \left\{ q_{\beta'}
        \Phi^{[j]}  +\Phi^{[j+1]}   \right\} (u_n)   \right) \\
        &  \leq  & \sum_{j=1}^{k-1} \pr_{f,s}  \left(  |\hat  \Phi^p_n  (u_n)  -\Phi^p(u_n)  |  \geq
          \frac{1}{2}  \left\{ q_{\beta'}  \Phi^{[j]} +  \Phi^{[j+1]} \right\}  (u_n) -
          |\Phi^p(u_n)| \right) .
\end{eqnarray*}
As $|\Phi^p(u_n)| \geq q_{\beta'}(u_n) \Phi^g (u_n)$ for large enough
$n$, we get
  \begin{eqnarray*}
   \pr_{f,s}(\hat s_n < s_k)
&\leq&  \sum_{j=1}^{k-1} \pr_{f,s}  \left(|\hat  \Phi^p_n   (u_n)  -\Phi^p(u_n)   |\geq  \frac{1}{2}
 \left\{ q_{\beta'}  \Phi^{[j]} +  \Phi^{[j+1]} \right\}  (u_n)  (1+o(1))\right) \\
  &\leq &  \sum_{j=1}^{k-1} \exp \left[ -\frac{n}{4}\left( A^2 u_n^{-2\beta'} \exp(-2u_n^{s_j} )+ \exp(-2 u_n^{s_{j+1}}) \right) \right]\\
  &\leq & N \exp \left( -\frac{n}{4}  \exp(-2u_n^{s}) \right).
\end{eqnarray*}
Now consider  the case  $\hat s_n >s_k $.
\begin{multline*}
   \pr_{f,s}(\hat  s_n  >  s_k)    \leq   \sum_{j=k+1}^{N}\pr_{f,s}  \left(  |\hat  \Phi_n^p(u_n)|  \leq
    \frac{1}{2} \{ q_{\beta'} \Phi^{[j-1]} + \Phi^{[j]} \} (u_n) \right) \\
     \leq \sum_{j=k+1}^{N} \pr_{f,s} \left( |\hat \Phi_n^p(u_n)- \Phi^p(u_n)|
     \geq  q_{\beta'}(u_n) \Phi^g (u_n) -
    \frac{1}{2} \{ q_{\beta'} \Phi^{[j-1]} + \Phi^{[j]} \} (u_n) \right) \\
    \leq N  \pr_{f,s} \left( |\hat \Phi_n^p(u_n)- \Phi^p(u_n)| \geq  q_{\beta'}(u_n) \{\Phi^g (u_n) -
    \frac{1}{2}  \Phi^{[k]}(u_n) \} +o(q_{\beta'}(u_n) \Phi^g(u_n))  \right)
  \end{multline*}
as $|\Phi^p(u_n)| \geq q_{\beta'}(u_n) \Phi^g (u_n)$ for large enough $n$  and $j-1\geq k$ .
According to the choice of the grid, we have $|s-s_k|\leq d_n $ and $d_n \log u_n \to 0$, which implies
\begin{eqnarray*}
\Phi^g(u_n)-\frac{1}{2} \Phi^k(u_n) &=& \exp(-u_n^s) \left( 1-\frac{1}{2} \exp[u_n^{s_k}(u_n^{s-s_k}-1 )] \right) \\
&=& \exp(-u_n^s) \left( 1-\frac{1}{2} \exp [u_n^{s_k}(s-s_k) \log u_n (1+o(1)) ] \right) \\
&\geq & \exp(-u_n^s) \left( 1-\frac{1}{2} \exp (u_n^{s_k-\smax} (1+o(1)) ) \right) \\
&\geq & \frac{1}{2} \exp(-u_n^s) (1+o(1) )  ,
 \end{eqnarray*}
 where the first inequality comes from $d_n \log u_n \leq u_n^{-\smax}$.
This gives
\begin{eqnarray*}
   \pr_{f,s}(\hat  s_n  >  s_k)
   & \leq  & N  \pr_{f,s} \left( |\hat \Phi_n^p(u_n)- \Phi^p(u_n)| \geq    \frac{1}{2} q_{\beta'}(u_n) \Phi^g (u_n)(1+o(1) ) \right) \\
   & \leq & N \exp \left( -\frac{A^2}{2} n u_n^{-2\beta'} \exp(-2u_n^s) (1+o(1) )\right).
  \end{eqnarray*}
In conclusion, as soon as we have $d_n \log u_n \leq u_n^{-\smax}$, and
$\log N = o( (\log n)^\alpha)$ (which is ensured by our choice of $d_n$) we get, for any $\epsilon>0$ and large enough $n$,
\begin{multline*}
   \pr_{f,s}(|\hat  s_n  -  s| \geq\epsilon)    \leq   N \exp \left( -\frac{A^2}{2} n u_n^{-2\beta'} \exp(-2u_n^{s}) (1+o(1) )\right) \\
    \leq  \exp \left( -\frac{A^2}{2}(\log n)^{\alpha}(1+o(1) ) \right).
  \end{multline*}
The last term gives a convergent series and then according to Borel
Cantelli's lemma, $\pr_{f,s}(|\hat  s_n  -  s| \geq\epsilon \text{
i.o }) =0$ leading to the almost sure convergence of $\hat s_n$.
\end{proof}
\\

\begin{proof}[Proof of Corollary~\ref{cor:deconv_s_inconnu}]
Note that the new choice of $d_n$ still satisfies the requirements for Proposition~\ref{conv_s} to be valid.
We introduce  respectively, $  h_n$, the non-random  version of  the bandwidth
$\hat h_n$  and $K_n$  the non-random  version of the  kernel $\hat  K_n$ both
constructed  with  self-similarity   index  $\sgrid$.  The  Fourier  transform
$\Phi^{K_n} $ of $K_n$ thus satisfies 
\begin{eqnarray*}
 \Phi^{K_n}(u) &=& \exp ((|u|/h_n)^{s_n(s)}) 1_{|u| \leq 1} \\
\text{where }  h_n &=& ( 2 ^{-1}\log  n -(\betasup-\sgrid +1/2) \log  \log n /
\sgrid )^{-1/ \sgrid}.  
\end{eqnarray*}
We also introduce the  corresponding  (classical) estimator
$$ f_n (x) = (n  h_n)^{-1} \sum_{i=1}^n K_n( h_n^{-1}(x-Y_i)).$$
Note that obviously, $\sgrid, K_n$ and $ h_n$ are unknown to the statistician.
These objects are used only as tools to assess the convergence of the procedure.
Now, remark that we have
\begin{multline*}
\esp_{f,s} [| \hat f_n(x) -f(x)|^2] = \esp_{f,s}  [|  f_n(x) -f(x)|^2 1_{\hat s_n =\sgrid}] +
\esp_{f,s} [| \hat f_n(x) -f(x)|^2 1_{\hat s_n \neq \sgrid}] \\
=T_1 +T_2,
\end{multline*}
say. Let us focus on the first term
$$ T_1 \leq \esp_{f,s} [|  f_n(x) -f(x)|^2 ] =\{\esp_{f,s}[ f_n(x)] -f(x)\}^2 + \var_s \{ f_n(x) \},
$$
introducing the bias and the variance of the estimator $ f_n(x)$. The important thing to note is that the kernel estimator $f_n$ uses parameter $\sgrid$ which is not equal to the true one $s$. Thus $T_1$ is not the classical risk for kernel estimator with known index $s$. Using Parseval's equality
\begin{multline*}
\{\esp_{f,s}[ f_n(x)] -f(x)\}^2
= \frac{1}{4\pi} \left[
\int e^{-iux} \Phi(u) \left( 1_{|u|\leq 1/ h_n} \exp(-|u|^s+|u|^{\sgrid}) -1 \right) du \right]^2\\
\leq  \frac{1}{4\pi} \left[ \int_{|u|\leq 1/ h_n} |\Phi(u)| \left(\exp(-|u|^s+|u|^{\sgrid})-1 \right) du
 +  \int_{|u|> 1/ h_n} |\Phi(u)| du  \right]^2.
\end{multline*}
The second term in the right hand side is the classical bias and equals $O( h_n^{\beta-1/2})$.
As soon as $d_n  h_n^{-s} \log  (1/ h_n) $ converges to zero, we can use the following development in the first term, uniformly for $|u| \leq 1/h_n$,
\begin{eqnarray*}
\exp(-|u|^s+|u|^{\sgrid})-1 &=& \exp\{ |u|^s(\sgrid-s)\log |u| (1+o(1)) \} -1 \\
&=& |u|^s (\sgrid-s)\log |u| (1+o(1)),
\end{eqnarray*}
which leads to
\begin{eqnarray*}
&& \{\esp_{f,s} [f_n(x)] -f(x)\}^2 \\
&\leq & \frac{1}{4\pi}
\left[ \left(  \int_{|u|\leq 1/  h_n} |\Phi(u)| |u|^s (\sgrid-s)\log |u| du \right) (1+o(1))
+  O(  h_n^{\beta-1/2})  \right]^2 \\
&\leq & O(d_n^2) 1_{\beta > s+1/2}
 + O(d_n^2   h_n^{2\beta-2s-1} \log^2(1/  h_n))1_{\beta \leq s+1/2} + O(  h_n^{2\beta-1}) .
\end{eqnarray*}
It can be easily seen that
\begin{eqnarray*}
 && d_n ^2   h_n^{-2s} \log^2(1/  h_n)\\
 &\leq& \frac {O(1) (\log n)^{2s/ \sgrid} (\log \log n)^2} { \log^2 n (\log \log n)^2 }
 =  O(1) (\log n)^{2(s-\sgrid)/\sgrid} \\
 &\leq& O(1) (\log n)^{2d_n/ \smin} \leq O(1) \exp\{2 \smax /(\smin \log n)\} =O(1),
\end{eqnarray*}
leading to
\begin{equation*}
\{\esp_{f,s}[ f_n(x)] -f(x)\}^2 \leq O(d_n^2) 1_{\beta > s+1/2}
+ O(  h_n^{2\beta-1}) .
\end{equation*}
Moreover, when $\beta > s+1/2$, we use $d_n^2 \leq (\log n)^{-(2\betasup -1)/ \smin} =O(  h_n^{(2\beta -1)} )$. With this choice of $d_n$, we thus ensure that in any case
$$ \{\esp_{f,s}[ f_n(x)] -f(x)\}^2 \leq O(  h_n^{2\beta-1}) .
$$
The variance of $ f_n(x)$ is bounded by
\begin{multline*}
 \var_{f,s} \{  f_n(x) \} = \frac{1}{4\pi^2 n } \esp_{f,s} \left[  \int_{|u| \leq 1/  h_n} e^{-iux} e^{|u|^{\sgrid}} (e^{iuY} -\Phi^p(u)) du \right]^2 \\
\leq \frac{1}{\pi^2 n }  \left( \int_{|u| \leq 1/ \bar h_n} e^{|u|^{\sgrid}} du \right)^2
= O\left( \frac{   h_n^{2(\sgrid -1)} \exp(2/  h_n^{\sgrid})} {n }\right) .
\end{multline*}
We finally get the bound
$$
T_1 \leq   O(  h_n^{2 \beta -1}) + O\left( \frac{  h_n^{2(\sgrid -1)} \exp(2/  h_n^{\sgrid})} {n }\right) .
$$
Now, we  prove that  the second term $T_2$ is negligible in front of the main term $T_1$, by  using Proposition~\ref{conv_s} and uniform bounds  on $| \hat
f_n (x)|$ and $|f(x)|$. First,
\begin{eqnarray*}
| \hat f_n (x)|
& \leq &\int e^{|t|^{\smax}} 1_{|t|\leq 1/\hat h_n} dt = O(\hat h_n^{\smax-1} \exp\{1/\hat h_n^{\smax}\}) \\
&\leq & O(1) (\log n)^{(1-\smax)/\smin} \exp\{ (\log n )^{\smax /\smin} \}   \\
| f (x)|  & \leq  &\int |\Phi(t)| dt  = O(\int  (1+|t|^{ 2\beta} )  ^{-1} dt)
=O(1),
\end{eqnarray*}
and then
\begin{eqnarray*}
T_2  &=& O(  (\log n)^{2(1-\smax)/\smin} \exp\{ 2(\log n )^{\smax /\smin} \} )  \pr_{f,s} (\hat   s_n  \neq  \sgrid)  \\
&= &  O\left( (\log n)^{2(1-\smax)/\smin}    \exp\left(2 (\log n )^{\smax /\smin} -\frac{A^2}{4} (\log n)^\alpha  (1+o(1))\right)   \right).
\end{eqnarray*}
As soon as we choose $\alpha >\smax /\smin$, this second term $T_2$ will be negligible in front of $T_1$.
In conclusion,
\begin{eqnarray*}
\esp_{f,s} [| \hat f_n(x) -f(x)|^2] &=& O(  h_n^{2 \beta -1}) + O\left(   h_n^{2(\sgrid -1)} \frac{\exp(2/  h_n^{\sgrid})} {n }\right)  \\
&=& O( (\log n)^{-(2\beta -1)/\sgrid}) = O( (\log n)^{-(2\beta -1)/s}) .
\end{eqnarray*}
\end{proof}

\bibliographystyle{natbib}
\bibliography{ccc}

\begin{thebibliography}{}

\bibitem[Butucea(2004a)Butucea]{ButuceaSORT}
Butucea, C. (2004a).
\newblock {Asymptotic normality of the integrated square error of a density
  estimator in the convolution model}.
\newblock {\em SORT\/}, {\bf 28}(1), 9--26.

\bibitem[Butucea(2004b)Butucea]{Butucea-deconv}
Butucea, C. (2004b).
\newblock {Deconvolution of supersmooth densities with smooth noise}.
\newblock {\em Can. J. Stat.}, {\bf 32}(2), 181--192.

\bibitem[Butucea(2004c)Butucea]{Butucea}
Butucea, C. (2004c).
\newblock {Quadratic functional estimation in view of minimax goodness-of-fit
  testing from noisy data}.
\newblock {\em To appear in Ann. Stat.}
\newblock Available at \texttt{arXiv:math.ST/0612361} (short version) and
  \texttt{http://www.proba.jussieu.fr/pageperso/butucea/RevAOSlong.pdf} (long
  version).

\bibitem[Butucea and Matias(2005)Butucea and Matias]{Butmat}
Butucea, C. and Matias, C. (2005).
\newblock Minimax estimation of the noise level and of the deconvolution
  density in a semiparametric convolution model.
\newblock {\em Bernoulli\/}, {\bf 11}(2), 309--340.

\bibitem[Butucea and Tsybakov(2007)Butucea and Tsybakov]{ButTsy}
Butucea, C. and Tsybakov, A.~B. (2007).
\newblock {Sharp optimality for density deconvolution with dominating bias. I.
  II}.
\newblock {\em Theory Probab. Appl\/}, {\bf 52}.
\newblock To appear.

\bibitem[Carroll and Hall(1988)Carroll and Hall]{Carroll-Hall}
Carroll, R.~J. and Hall, P. (1988).
\newblock Optimal rates of convergence for deconvolving a density.
\newblock {\em J. Amer. Statist. Assoc.}, {\bf 83}(404), 1184--1186.

\bibitem[Comte {\em et~al.}(2006)Comte, Rozenholc, and
  Taupin]{Comte-Taupin-Rozen}
Comte, F., Rozenholc, Y., and Taupin, M.-L. (2006).
\newblock Penalized contrast estimator for adaptive density deconvolution.
\newblock {\em Can. J. Stat.}, {\bf 34}(3), 431--452.

\bibitem[Efromovich(1997)Efromovich]{Efromovich}
Efromovich, S. (1997).
\newblock Density estimation for the case of supersmooth measurement error.
\newblock {\em J. Amer. Statist. Assoc.}, {\bf 92}(438), 526--535.

\bibitem[Fan(1991)Fan]{Fan1}
Fan, J. (1991).
\newblock On the optimal rates of convergence for nonparametric deconvolution
  problems.
\newblock {\em Ann. Statist.}, {\bf 19}(3), 1257--1272.

\bibitem[Fan(1992)Fan]{Fan3}
Fan, J. (1992).
\newblock {Deconvolution with supersmooth distributions}.
\newblock {\em Can. J. Stat.}, {\bf 20}(2), 155--169.

\bibitem[Fan and Koo(2002)Fan and Koo]{Fan-Koo}
Fan, J. and Koo, J.-Y. (2002).
\newblock {Wavelet deconvolution}.
\newblock {\em IEEE Trans. Inf. Theory\/}, {\bf 48}(3), 734--747.

\bibitem[Gayraud and Pouet(2005)Gayraud and Pouet]{GayraudPouet}
Gayraud, G. and Pouet, C. (2005).
\newblock Adaptive minimax testing in the discrete regression scheme.
\newblock {\em Probab. Theory Relat. Fields\/}, {\bf 133}(4), 531--558.

\bibitem[Hall(1984)Hall]{Hall}
Hall, P. (1984).
\newblock {Central limit theorem for integrated square error of multivariate
  nonparametric density estimators}.
\newblock {\em J. Multivar. Anal.}, {\bf 14}(3), 1--16.

\bibitem[Hall and Heyde(1980)Hall and Heyde]{HallHeyde}
Hall, P. and Heyde, C. (1980).
\newblock {\em {Martingale limit theory and its application}\/}.
\newblock {Probability and Mathematical Statistics. Academic Press}, {New
  York}.

\bibitem[Holzmann {\em et~al.}(2007)Holzmann, Bissantz, and Munk]{HolzBisMunk}
Holzmann, H., Bissantz, N., and Munk, A. (2007).
\newblock {Density testing in a contaminated sample}.
\newblock {\em J. Multivariate Analysis\/}, {\bf 98}(1), 57--75.

\bibitem[Ingster and Suslina(2003)Ingster and Suslina]{Ingster}
Ingster, Y. and Suslina, I. (2003).
\newblock {\em {Nonparametric goodness-of-fit testing under Gaussian
  models}\/}.
\newblock {Lecture Notes in Statistics. 169. Springer.}, New York.

\bibitem[Korolyuk and Borovskikh(1994)Korolyuk and Borovskikh]{BorKor}
Korolyuk, V. and Borovskikh, Y. (1994).
\newblock {\em {Theory of $U$-statistics}\/}.
\newblock {Mathematics and its Applications. Kluwer Academic Publishers},
  {Dordrecht}.

\bibitem[Lehmann and Romano(2005)Lehmann and Romano]{LehRom}
Lehmann, E.~L. and Romano, J.~P. (2005).
\newblock {\em Testing statistical hypotheses\/}.
\newblock Springer Texts in Statistics. Springer, New York, third edition.

\bibitem[Matias(2002)Matias]{Matias}
Matias, C. (2002).
\newblock Semiparametric deconvolution with unknown noise variance.
\newblock {\em ESAIM, Probab. Stat.}, {\bf 6}, 271--292.
\newblock (electronic).

\bibitem[Meister(2004)Meister]{Meister}
Meister, A. (2004).
\newblock {On the effect of misspecifying the error density in a deconvolution
  problem}.
\newblock {\em Can. J. Stat.}, {\bf 32}(4), 439--449.

\bibitem[Pensky and Vidakovic(1999)Pensky and Vidakovic]{penvid}
Pensky, M. and Vidakovic, B. (1999).
\newblock Adaptive wavelet estimator for nonparametric density deconvolution.
\newblock {\em Ann. Statist.}, {\bf 27}(6), 2033--2053.

\end{thebibliography}

\newpage
%
%
\appendix

\section{Technical Proofs}

\begin{proof} [Proof of Lemma~\ref{lem2}]
Using a Markov inequality and the usual controls on bias and variance, we get
\begin{multline*}
 \pr_0(|T_{n,N+1} -\esp_0(T_{n,N+1})|  >  \mathcal{C}^\star
t_{n,N+1}^2 -\esp_0(T_{n,N+1})) \leq \frac{C n^{-2} (h^{N+1})^{-(4\sigma +1)}}
{(\mathcal{C}^{\star}    t_{n,N+1}^2   -   C    (h^{N+1})^{2\betasup}   )^{2}}
\\
=O(\frac{1}{\mathcal{C}^\star-C}),
\end{multline*}
and  by choosing  $\mathcal{C}^\star$ large  enough, this  term is
smaller than some $\epsilon >0$.
\end{proof}
\\

\begin{proof} [Proof of Lemma~\ref{lem3}]
Let us write
$$ \pr_f ( | T_{n,N+1}| \leq \mathcal{C}^\star t_{n,N+1}^2) \leq \pr_f (  | T_{n,N+1} -\esp_f T_{n,N+1} | \geq \|f-f_0\|_2^2 -\mathcal{C}^\star t_{n,N+1}^2-B_f(T_{n,N+1}))
$$
where
\begin{eqnarray*}
&& |B_f(T_{n,N+1})| = |\esp_f (T_{n,N+1}) -\|f-f_0\|_2^2| \\
&\leq &
\int_{|u|\geq  1/h^{N+1}} |\Phi(u)|^2  du +  2\left(  \int_{|u|\geq 1/h^{N+1}}
  |\Phi(u)|^2 du \int_{|u|\geq 1/h^{N+1}} |\Phi_0(u)|^2 du \right)^{1/2}\\
& \leq& \big(L (h^{N+1})^{2\beta} \exp\{ -2\alpha /(h^{N+1}) ^r\} +2L(h^{N+1})^{\beta+\betasup} \exp\{ -\alpha /(h^{N+1}) ^r\}
\big)(1+o(1))\\
&\leq & 2L(h^{N+1})^{\beta+\betasup} \exp\{ -\alpha /(h^{N+1}) ^r\} (1+o(1)).
\end{eqnarray*}
In the same way as in the proof of Lemma~\ref{lem4}, we have
$$
 \esp_f (T_{n,N+1} -\esp_f(T_{n,N+1}))^2
\leq \frac{C }{n^2 (h^{N+1})^{4\sigma +1}} +\frac{4 \Omega_g^2(f-f_0)} {n} 1_{\beta\geq \sigma}
= w_{n,f}^2,
$$
and $\Omega_g(f-f_0)$  is a constant  depending on $f$  and $g$ (but not  $n$) and
satisfying \linebreak $|\Omega_g^2(f-f_0) | \leq C\|f- f_0\|_2^{2-2\sigma /\betasup}$. 
The rest of the proof follows the same lines as Lemma~\ref{lem4}.
Indeed, Markov's Inequality leads the following bound on the second type error term
\begin{multline*}
\frac{w_{n,f}^2}{(\|f-f_0\|_2^2         -\mathcal{C}^\star         t_{n,N+1}^2
  -2L(h^{N+1})^{2\beta} \exp\{-\alpha /(h^{N+1})^r\} (1+o(1)))^2}\\
\leq     \max\left(      \frac{Cn^{-2}(h^{N+1})^{-4\sigma   -1}}     {(\mathcal{C}^0
  -\mathcal{C}^\star   )^2  \psi_{n,r}^{4}}   ;  \frac{C} {n \|f-f_0\|_2^{2+2\sigma
    /\betasup }(\mathcal{C}^0-\mathcal{C}^\star )^2 }
\right)
\end{multline*}
The first term in  the right hand side is a constant which  can be as small as
we need, by choosing a  large enough constant $\mathcal{C}^0$. The second term
converges to zero.
\end{proof}
\\

\begin{proof}[Proof of Lemma~\ref{lem5}]
As $\beta > \nu$,  the bandwidths satisfy $h_\nu h_\beta^{-1}
=o(1)$.  Then, as $G$ is compactly supported on $[-1,0]$, we have
\begin{multline*}
\esp_{0} \left(
  \frac {G_{\beta}\left( Y_1  -x_{j,\beta} \right)G_{\nu}\left( Y_1 -x_{i,\nu}
    \right)} {p_0^2\left( Y_1 \right)} \right)   =
\int_{\mathbb R} \frac {G_{\beta}\left( y -x_{j,\beta} \right)
G_{\nu}\left( y -x_{i,\nu} \right)} {p_0\left( y \right)}
\; dy \\
 =  \int_{\left[-1, \; 0 \right]}
\frac {G_{\beta}\left( h_\nu u + x_{i,\nu} -x_{j,\beta} \right)
G\left( u \right)} {p_0\left(  h_\nu u + x_{i,\nu} \right)}
\; du.
\end{multline*}
Apply the Taylor Formula to get
\begin{eqnarray*}
G_{\beta}\left( h_\nu u + x_{i,\nu} -x_{j,\beta} \right)
 & = & G_{\beta}\left(  x_{i,\nu} -x_{j,\beta} \right)
+ \frac { h_\nu } {h_\beta^2} u
G'\left( \frac {h_\nu {\tilde u_1} + x_{i,\nu} -x_{j,\beta}}
{h_\beta}  \right) \\
\text{and } \; \; \; \frac 1 {p_0\left( h_\nu u + x_{i,\nu} \right)}
& = & \frac 1 {p_0 \left(  x_{i,\nu} \right)}
- \frac { p_0'
\left(  h_\nu {\tilde u_2} + x_{i,\nu} \right)}
{p_0 \left(  h_\nu {\tilde u_2} + x_{i,\nu} \right)^2 }
h_\nu u ,
\end{eqnarray*}
where $0 \leq \tilde u_1 \leq u$ and $0 \leq \tilde u_2 \leq u$.
As $\int G=0$, we obtain
\begin{eqnarray*}
&& \int_{\left[-1, 0 \right]}
\frac {G_{\beta}\left( h_\nu u + x_{i,\nu} -x_{j,\beta} \right)
G\left( u \right)} {p_0\left(  h_\nu u + x_{i,\nu} \right)}
\; du\\
& =& \frac 1 {p_0 \left( x_{i,\nu} \right)} \frac { h_\nu }
{h_\beta^2} \int_{\left[-1,  0 \right]} u G'
\left( \frac {h_\nu {\tilde u_1} + x_{i,\nu} -x_{j,\beta}}
{h_\beta}  \right) G(u) du \\
&& - h_\nu G_{\beta}\left( x_{i,\nu} -x_{j,\beta} \right)
\int_{\left[-1,  0 \right]}
\frac { p_0'
\left(  h_\nu {\tilde u_2} + x_{i,\nu} \right)}
{p_0 \left(  h_\nu {\tilde u_2} + x_{i,\nu} \right)^2 }
u G\left( u \right)  du \\
&& \mathop{-} \frac { h_\nu^2 } {h_\beta^2} \int_{\left[-1, 0 \right]}
\frac { p_0'
\left(  h_\nu {\tilde u_2} + x_{i,\nu} \right)}
{p_0 \left(  h_\nu {\tilde u_2} + x_{i,\nu} \right)^2 }
 u^2 G'
\left( \frac {h_\nu {\tilde u_1} + x_{i,\nu} -x_{j,\beta}}
{h_\beta} \right) G\left( u \right) \; du .
\end{eqnarray*}
This leads to
$$
 \esp_{0} \left(
  \frac {G_{\beta}\left( Y_1  -x_{j,\beta} \right)G_{\nu}\left( Y_1 -x_{i,\nu}
    \right)} {p_0^2\left( Y_1 \right)}
\right) = \frac{h_\nu}{h_\beta^2} R_{i,j}
$$
where
\begin{eqnarray*}
R_{i,j} &=& \frac 1 {p_0  \left( x_{i,\nu} \right)} \int_{\left[-1, 0 \right]} u
G' \left( \frac {h_\nu {\tilde u_1} + x_{i,\nu} -x_{j,\beta}}
{h_\beta} \right) G(u) du \\
&& - h_\beta G\left(\frac{ x_{i,\nu}
-x_{j,\beta}}{h_\beta} \right) \int_{\left[-1,  0 \right]}
\frac { p_0'
\left(  h_\nu {\tilde u_2} + x_{i,\nu} \right)}
{p_0 \left(  h_\nu {\tilde u_2} + x_{i,\nu} \right)^2 }
u G\left( u \right) du \\
&& - h_\nu \int_{\left[-1, 0 \right]}
\frac { p_0'
\left(  h_\nu {\tilde u_2} + x_{i,\nu} \right)}
{p_0 \left(  h_\nu {\tilde u_2} + x_{i,\nu} \right)^2 }
 u^2 G'
\left( \frac {h_\nu {\tilde u_1} + x_{i,\nu} -x_{j,\beta}}
{h_\beta} \right) G\left( u \right) \; du .
\end{eqnarray*}
satisfies
$$
|R_{ij}|  \leq  (\inf_{[0,1]}  p_0)^{-1} \|G\|_\infty  \|G'\|_\infty  +
\|G\|_\infty  \|p_0'\|_\infty  (\inf_{[-1,1]}  p_0)^{-2} (h_\beta  \|G\|_\infty  +h_\nu
\|G'\|_\infty) ,
$$
which ends the proof of Lemma~\ref{lem5}.
\end{proof}
\\


\begin{proof} [Proof of Theorem~\ref{th:bruit_poly_f0ana}]

Assume now that $f_0 \in \mathcal{F}(\alphamax, \rsup, \beta_0, L)$,
for some $\beta_0 \in [\betainf, \betasup]$. The proof follows the
same lines as the proof of Theorem~\ref{th:bruit_poly_f0sob}.

For the first-type error we write
\begin{eqnarray*}
  \pr_0(\Delta_n^*=1) &=& \sum_{i=1}^{N_1}
  \pr_0(|T_{n,i} -\esp_0(T_{n,i})|  >  \mathcal{C}^\star t_{n,i}^2 -\esp_0(T_{n,i}))
\\
  && + \sum_{i-N_1=1}^{N_2} \pr_0(|T_{n,i} -\esp_0(T_{n,i})| > \mathcal{C}^\star t_{n,i}^2 -\esp_0(T_{n,i})).
\end{eqnarray*}
For the first $N_1$ terms we apply Lemma~\ref{lem1} with
$\esp_0(T_{n,i}) = o(1) L (h_i)^{2 \beta_0} \exp(-2
\alphamax/h_i^{\rsup})$ which is smaller than $t_{n,i}^2$ for all
$i=1,\dots,N_1$ and the same result follows. For the last $N_2$
terms we also use the Berry-Esseen inequality as in the proof of
Lemma~\ref{lem1} for
$$
x=\mathcal{C}^\star t_{n,i}^2 -\esp_0(T_{n,i})\geq \mathcal{C}^\star t_{n,i}^2 (1-o(1))
$$
as $\esp_0(T_{n,i}) = o(1) h_i^{2 \beta_0}\exp(-2\alphamax/h_i^{\rsup}) = o(1/n)$. We get
$x/v_n =O(1) (\log \log \log n)^{1/2}$
\begin{eqnarray*}
&&\sum_{i-N_1=1}^{N_2} \pr_0(|T_{n,i} -\esp_0(T_{n,i})| > \mathcal{C}^\star t_{n,i}^2 -\esp_0(T_{n,i}))\\
&\leq & N_2 \frac{v_n}{\mathcal{C}^\star t_{n,i}^2}
\exp\left( - \frac{(\mathcal{C}^\star)^2 t_{n,i}^4 }{4 v_n^2}\right)
\leq C_1 \frac{(\log \log \log n)^{-1/2} }{(\log \log n)^{b-1}}=o(1),
\end{eqnarray*}
for some $b>1$ for $\mathcal{C}^\star$ large enough.
Indeed, all other calculations are similar as they are related mostly to the
distribution of the noise which didn't change.

As for the second-type error,
\begin{eqnarray*}
&& \sup_{\tau \in \mathcal{T}} \sup_{f \in \mathcal{F}(\tau,L)} \pr_f(\Delta_n^\star =0)\\
& \leq & 1_{\rinf =\rsup=0} \sup_{\alpha\geq \alphamin, \beta \in [\betainf;\betasup]}
\sup_{\substack{f \in  \mathcal{F}(\alpha,0,\beta,L)\\ \|f-f_0\|_2^2 \geq
 \mathcal{C} \psi^2_{n,(\alpha,0,\beta)}}}
\pr_f ( \forall 1\leq i \leq N_1, \; | T_{n,i}| \leq \mathcal{C}^\star t_{n,i}^2) \\
& + & 1_{\rinf >0} \sup_{ r \in [\rinf;\rsup], \alpha\in [\alphamin,\alphamax], \beta\in [\betainf,\betasup]}
\sup_{\substack{f \in  \mathcal{F}(\tau,L)\\ \|f-f_0\|_2^2 \geq
 \mathcal{C} \psi^2_{n,\tau}}}
\pr_f ( \forall N_1+1\leq i \leq N_1+N_2, \; | T_{n,i}| \leq \mathcal{C}^\star t_{n,i}^2).
\end{eqnarray*}
For the first term in the previous sum we actually apply precisely
Lemma~\ref{lem4}. For the second term we mimic the proof of
Lemma~\ref{lem4} and choose some $f$ in $\mathcal{F}(\alpha, r,
\beta,L)$ such that $\|f-f_0\|_2^2 \geq \mathcal{C} \psi^2_{n,r}$,
where we denote $\psi_{n,r}=\psi_{n,\tau} 1_{r>0}$. We define $r_f$
as the smallest point on the grid $\{r_1,\ldots, r_{N_2}\}$ such
that $r \leq r_f$. We denote by $h_f$, $t^2_{n,f}$ and $T_{n,f}$ the
bandwidth, the threshold and the test statistic associated to
parameters $\alphamax$ and $r_f$ (they do not depend on $\beta$).
Then
\begin{eqnarray}
&&\pr_f ( \forall N_1+1\leq i \leq N_1+N_2, \; | T_{n,i}| \leq \mathcal{C}^\star t_{n,i}^2)
\nonumber \\
&\leq & \pr_f(|T_{n,f}-\esp_f(T_{n,f})|\geq \|f-f_0\|_2^2 -
\mathcal{C}^\star t^2_{n,f}- B_f(T_{n,f})), \label{prem}
\end{eqnarray}
where, as in Theorem~\ref{th:bruit_poly_f0sob}
\begin{eqnarray*}
|B_f(T_{n,f})| &=& |\|J_h * f - f\|_2^2 + 2 \langle f -J_h* f, f_0 \rangle |\\
&\leq & \big( L  h_f^{2 \beta }\exp(-2 \alpha/ h_f^r)
+2Lh_f^{\beta+\beta_0}\exp(-\alpha/h_f^r -\alphamax/h_f^{\rsup})
\big) (1+o(1)) \\
&\leq &L (h_f^{2 \beta }+h_f^{\beta+\beta_0})\exp(-2 \alpha/ h_f^r)(1+o(1))\\
&\leq &L (h_f^{ \beta +\beta\wedge \beta_0})\exp(-2 \alpha/ h_f^r)(1+o(1)).
\end{eqnarray*}
Using Markov's inequality, we get the following upper bound for
(\ref{prem})
\begin{equation}
\frac{\var_f(T_{n,f})}{(\|f-f_0\|_2^2 - \mathcal{C}^\star t^2_{n,f}-
B_f(T_{n,f}))^2}. \label{deuze}
\end{equation}
The variance is bounded from above by
\begin{equation}\label{var}
\esp_f(T_{n,f}-\esp_f(T_{n,f}))^2\leq \frac{C}{n^2 h_f^{4\sigma +1}}+ \frac{4 \Omega_g^2(f-f_0)}{n},
\end{equation}
and similarly to \citep{Butucea}     we      show     that
$\Omega_g^2(f-f_0)     \leq \|f-f_0\|_2^2 (\log
\|f-f_0\|_2^{-2})^{2\sigma/r}$.
We have
\begin{equation*}
t_{n,f}^2 \psi_{n,r}^{-2} =(\log n)^{(4\sigma+1)(1/r_f-1/r)/2} \leq
1,
\end{equation*}
and thus $\|f-f_0\|_2^2 -\mathcal{C}^\star t_{n,f}^2 \geq
(\mathcal{C} -\mathcal{C}^\star) \psi_{n,r}^2$. Moreover,
\begin{multline*}
B_f(T_{n,f}) \psi_{n,r}^{-2} \leq C (\log \log \log n)^{-1/2} \left(
\log n \right)^{-(\beta+\beta \wedge\beta_0)/r_f
-(4\sigma+1)/(2r)}\\
\times \exp\left\{-2\alpha \left(\frac{\log n}{2c} \right)^{r/r_f}
+\log n \right\}.
\end{multline*}
The construction of the grid ensures that $-1/(\log \log n) \leq
r-r_f \leq 0$ and thus
\begin{eqnarray*}
&&\exp\left\{-2\alpha \left(\frac{\log n}{2c} \right)^{r/r_f} +\log
n \right\} \\
&=& \exp\left\{- \frac{\log n}{c} \left[ \alpha
 \exp \left(\frac{r-r_f}{r_f}\log \log n (1+o(1)) \right) -c \right]
 \right\}\\
 &\leq& \exp\left\{- \frac{\log n}{c} \left[ \alphamin
 \exp \left(\frac{-1}{\rinf} (1+o(1)) \right) -c \right]
 \right\}=O(1),
\end{eqnarray*}
as we chose the constant $c < \alphamin \exp(-1/\rinf)$. Finally,
we have $B_f(T_{n,f})\psi_{n,r}^{-2}=o(1)$. Let us come back to
\eqref{deuze}. We distinguish two cases whether the first or the
second term in \eqref{var} is dominant. If the first term in the
variance dominates, we have the following bound for \eqref{deuze}
\begin{equation*}
    \frac{n^{-2}h_f^{-(4 \sigma+1)}}{(\mathcal{C}-\mathcal{C}^\star)^2 \psi_{n,\tau}^4}
    \leq \frac{C} {\log \log \log n} \to 0.
\end{equation*}
On the other  hand, if the second  term in (\ref{var}) is the larger
one, the bound \eqref{deuze} writes
\begin{multline*}
\frac{n^{-1}\|f-f_0\|_2^2 (\log \|f-f_0\|_2^{-2})^{2\sigma/r}}
{\|f-f_0\|_2^4 (1- \mathcal{C}^\star / \mathcal{C} +o(1))^2}  \leq C
n^{-1}\psi_{n,r}^{-2} (\log \psi_{n,r}^{-2})^{2\sigma/r} \\=C (\log
n)^{-1/(2r)}(\log \log \log n)^{-1/2} =o(1).
\end{multline*}
This finishes the proof.
\end{proof}
\\


\begin{proof}[Proof of Corollary~\ref{cor:intf2_s_inconnu}]
We keep on with the same notation as in Subsection~\ref{sec:apres_estim_s}
and denote by $I$ the functional $\int f^2$.
In the same way as in the proof of Corollary~\ref{cor:deconv_s_inconnu}, we write
\begin{equation} \label{toto}
 \esp_{f,s} [| \hat T_{n } -I |^2 ] \leq  \esp_{f,s} [| T_{n } -I  |^2]
 +\esp_{f,s} [| \hat T_{n } -I |^2 1 _{\{\hat s_n \neq \sgrid\}} ] .
\end{equation}
Let us first focus on the first term appearing in the right hand side of  \eqref{toto}. We split it into the square of a bias term plus a variance term. The bias is bounded by
\begin{eqnarray*}
| \esp_{f,s}   T_{n } -I  | &\leq& \frac{1}{2\pi} \left( \int_{|u| > 1/h_n } |\Phi(u)|^2 du + \int_{|u| \leq 1/h_n} | \exp(2|u|^{\sgrid} -2|u|^s) -1| |\Phi(u)|^2 du \right) \\
&\leq & O(h_n^{-2\beta}) +  \int_{|u| \leq 1/h_n} 2|u|^s |\sgrid -s | \log |u| |\Phi(u)|^2 du \\
&\leq & O(h_n^{-2\beta}) +O(d_n) 1_{\{s \leq 2\beta \}}  + O(h_n^{2\beta -s} \log(1/h_n) d_n)1_{\{s > 2\beta\}}.
\end{eqnarray*}
Like   in  the   proof  of   Corollary  \ref{cor:deconv_s_inconnu},   we  have
$d_nh_n^{-s}  \log(1/h_n)   =o(1)$  and  thus  using  that   $d_n  \leq  (\log
n)^{-2\betasup / \smin} =O((\log n)^{-2\beta/s})$, we finally get
$$
| \esp_{f,s}   T_{n } -I  | \leq O((\log n)^{-2\beta/s}) .
$$
Concerning the variance term, we easily get
$$
 \var_{f,s}   ( T_{n } ) \leq \frac{C_1}{n^2}   h_n^{\sgrid-1}\exp(4/  h_n ^{\sgrid}) + \frac {C_2}   n   h_n^{2\beta +\sgrid-1} \exp( 2/  h_n ^{\sgrid}) ,
$$
where $C_1$ and $C_2$ are positive constants (we refer to \citep{Butucea}, Theorem 4 for more details).
Using the form of the bandwidth $  h_n$, we have
$$
 \esp_{f,s} |    T_{n } -I  |^2 = O\left( \frac {\log n}{2} \right)^{-4\beta/ s} .
$$
Let us now focus on the second term appearing in the right hand side of  \eqref{toto}. Denoting by $h_0= (\log n /2)^{-1/\smin}$, we have
\begin{equation*}
| \hat T_n | \leq  \frac 1 {2\pi } \int_{|u| \leq 1/h_0} \exp(2|u|^{\smax} ) du
=O(h_0^{\smax -1} \exp(2/h_0^{\smax})).
\end{equation*}
Moreover,
$$
I = \|f\|_2^2 = \frac 1 {2\pi} \|\Phi \|_2^2
$$
This leads to
\begin{multline*}
\esp_{f,s} [| \hat T_{n } -I |^2 1_{ \{\hat s_n \neq \sgrid\}} ]
\leq  C\left ( \frac {\log n} {2} \right) ^{(1-\smax)/\smin}
\exp\left\{ 2\left( \frac {\log n }{2} \right) ^{\smax/\smin} \right\} \pr_{f,s}( \hat s_n \neq \sgrid)\\
\leq  C\left ( \frac {\log n} {2} \right) ^{(1-\smax)/\smin}
\exp\left\{ 2\left( \frac {\log n }{2} \right) ^{\smax/\smin} \right\}
 \exp\left( -\frac{A^2}{4} (\log n)^a  (1+o(1))\right) ,
\end{multline*}
and this term is negligible in front of the first term appearing in the right
hand side of \eqref{toto} as soon as $a > \smax/\smin$. This leads to the result.
\end{proof}
\\

\begin{proof}[Proof of Corollary~\ref{cor:test_s_inconnu}]
We use the same notation as in Subsection~\ref{sec:apres_estim_s}.
Moreover, $ T_n^0$ is the  test statistic constructed with the
deterministic kernel $ K_n$ and the deterministic bandwidth $ h_n$;
and $ t_n^2$ is the threshold defined with the  parameter value
$\sgrid$ for the  self-similarity index.  The first type error of
the test is controlled by
\begin{equation*}
\pr_{f_0,s}(\Delta_n ^\star =1) = \pr_{f_0,s} (|\hat T_n^0 |  \hat t_n^{-2} >
\mathcal{C}^\star) \\
\leq   \pr_{f_0,s}   ( \hat  s_n \neq  \sgrid)  + \pr_{f_0,s}  (|  T_n^0  | t_n^{-2}  >
\mathcal{C}^\star ) .
\end{equation*}
The first term on the right hand side of this inequality converges to zero according to Proposition \ref{conv_s}.
Let  us   focus  on   the  second  term.  We have
\begin{equation*}
\pr_{f_0,s} (| T_n^0 | t_n^{-2} > \mathcal{C}^\star )
\leq \frac{1}{(\mathcal{C}^\star) ^2 t_n^{4}} \esp_{f_0,s} ( T_{n }^0 )^2
\leq \frac{1}{(\mathcal{C}^\star) ^2 t_n^{4}} \left\{
(\esp_{f_0,s} T_{n }^0 )^2 + \var_{f_0,s} T_n^0 \right\} .
\end{equation*}
It is easily seen that
\begin{eqnarray*}
\esp_{f_0,s}  T_{n  }^0 &=&  \frac{1}{2\pi}  \int_{|u|\leq 1/h_n}  |\Phi_0(u)|^2
|\exp(|u|^{\sgrid}   -|u|^s)   -1|^2   du   +\frac{1}{2\pi}   \int_{|u|>1/h_n}
|\Phi_0(u)|^2 du \\
& \leq & \frac{d_n^2}{2\pi}  \left(\int_{|u|\leq 1/h_n}  |\Phi_0(u)|^2
 |u|^{2s} \log ^2 |u| du \right) (1+o(1)) + O(h_n^{2\betasup} ) \\
& \leq & O(d_n^2)1_{\betasup >s}  + O(h_n^{2\betasup} ) = O(h_n^{2\betasup}).
\end{eqnarray*}
the  inequalities being valid as soon as $h_n^{-s}d_n \log(1/h_n)$
converges to zero.
Like in  the  proof  of  Theorem  4  in
\citep{Butucea}, we can show that
$$
\var_{f_0,s}(T_n^0) \leq O(1) \frac{h_n^{\sgrid -1}}{n^2} \exp(4/h_n^{\sgrid})
+ O(1) \frac{h_n^{2\betasup +\sgrid-1} }{n}\exp(2/h_n^{\sgrid}).
$$
Finally, we get
\begin{eqnarray*}
&& \pr_{f_0,s} (| T_n^0 | t_n^{-2} > \mathcal{C}^\star )\\
&\leq &\frac{1}{(\mathcal{C}^\star) ^2 t_n^{4}} \left\{
O(h_n^{4\betasup}) + O(1) \frac{h_n^{\sgrid -1}}{n^2} \exp(4/h_n^{\sgrid})
+ O(1) \frac{h_n^{2\betasup +\sgrid-1} }{n}\exp(2/h_n^{\sgrid})
\right\} \\
&\leq& \frac{O(1)}{\mathcal{C}^\star}.
\end{eqnarray*}

Choosing $\mathcal{C}^\star$  large enough achieves  the control of  the first
error term. We now turn to the second
error term. Under  hypothesis $H_1(\mathcal{C}, \Psi_n)$, there  exists some $\beta$
such   that   $f$   belongs    to   $\mathcal{F}(0,0,\beta,L)$   and   $\|f-f_0\|_2^2   \geq
\mathcal{C}\psi_{n,\beta}$. We write
\begin{equation*}
\pr_{f,s}(\Delta_n  ^\star =0)  =  \pr_{f,s}(|\hat T_n^0|  \hat t_n^{-2}  \leq
\mathcal{C}^\star) \leq  \pr_{f,s}(\hat s_n \neq \sgrid)  + \pr_{f,s}(| T_n^0|
t_n^{-2} \leq \mathcal{C}^\star) .
\end{equation*}
As already  seen, the  first term in  the right  hand side of  this inequality
converges  to zero,  so we  only deal  with the  second one.  We  define
$B_{f,s} (T_n^0) = \esp_{f,s} T_n^0 - \|f-f_0\|_2^2$. Thus
\begin{multline}
  \label{eq:zorro}
\pr_{f,s}(|  T_n^0| t_n^{-2}  \leq \mathcal{C}^\star)  \leq  \pr_{f,s}(| T_n^0
-\esp_{f,s}  T_n^0| \geq  \| f-f_0\|_2^2   -\mathcal{C}^\star
t_n^{2} +B_{f,s} (T_n^0)) \\
\leq \frac{\var_{f,s}(T_n^0)}{(\|f-f_0\|_2^2 -\mathcal{C}^\star t_n^2
+B_{f,s}(T_n^0))^2}.
\end{multline}
We compute this bias term  $B_{f,s}  (T_n^0) $.
\begin{eqnarray*}
  B_{f,s}  (T_n^0) &=&  \frac{1}{2\pi} \int  |\exp(|u|^{\sgrid}  -|u|^s) \Phi(u)
  1_{|u|\leq     1/h_n}     -\Phi_0(u)|^2     du     -\frac{1}{2\pi}     \int
  |\Phi(u)-\Phi_0(u)|^ 2 du \\
&\leq  &  \frac{1}{2\pi} \int_{|u|\leq  1/h_n}  |[\exp(|u|^{\sgrid} -|u|^s)  -1]
\Phi(u)|^2du +\frac{1}{2\pi} \int_{|u|>1/h_n} |\Phi(u)|^2 du \\
&\leq   &  \frac{d_n^2}{2\pi}(1+o(1))\int_{|u|\leq   1/h_n}  |u|^{2s}\log^2|u|
|\Phi(u)|^2 du + O(h_n^{2\beta})\\
&\leq& O(d_n^2)1_{\beta>s} +O(h_n^{2\beta}) = O(h_n^{2\beta}).
\end{eqnarray*}
In   fact,  there   exists   some  constant  $C_1>0$ depending only on $L$ and on
the noise distribution  such   that $B_{f,s}(T_n^0) \leq C_1 h_n^{2\beta}$.
Under hypothesis $H_1(\mathcal{C},\Psi_n)$, we also have $\|f-f_0\|_2^2 \geq
\mathcal{C}\psi_{n,\beta}^2$. Thus,
\begin{eqnarray*}
\|f-f_0\|_2^2 -\mathcal{C}^\star t_n^2 +B_{f,s}(T_n^0) & \geq &
\mathcal{C}\left(\frac{\log   n}{2}  \right)^{-2\beta/s}  -\mathcal{C}^{\star}
\left(\frac{\log  n}{2}\right)^{-2\betasup/\sgrid}  -C_1  \left(\frac{\log
      n}{2}\right) ^{-2\beta/\sgrid} \\
&\geq & a \left(\frac{\log n}{2} \right)^{-2\beta/s}.
\end{eqnarray*}
where    $a   =\mathcal{C}-\mathcal{C}^\star-C_1$    is    positive   whenever
$\mathcal{C}> \mathcal{C}^0 :=\mathcal{C}^\star -C_1$.
Returning to \eqref{eq:zorro}, we get
\begin{equation*}
  \pr_{f,s}(| T_n^0| t_n^{-2} ) \leq \frac{\psi_{n,\beta}^4}{a^2}
\var_{f,s}(T_n^0) .
\end{equation*}
Computation  of  the variance  follows  the  same  lines as  under  hypothesis
$H_0$.  We obtain
$$
\var_{f,s}(T_n^0)   \leq O(1) \frac{h_n^{\sgrid  -1}}{n}\exp(2/h_n^{\sgrid})\left(
  h_n^{2\beta} +\frac{\exp(2/h_n^{\sgrid})} {n} \right).
$$
The choice of the bandwidth ensures  that the second type error term converges
to zero.
\end{proof}
\\

\begin{proof}[Proof of Theorem~\ref{th:Ustats}]
  This  proof  follows  the  lines of  Theorem 3.9  in
\citep{HallHeyde}.
Combining   the   Skorokhod   representation   Theorem  and   Lemma   3.3   in
\citep{HallHeyde}, there exists a
nonnegative random  variable $T_n$ such that  for any $0<\epsilon<1/2$  and any real
$x$,
$$
|\pr(U_n  \leq x)  -\phi(x)|=  |\pr(S_n  \leq v_n^{-1}  x)  -\phi(x/v_n)| \leq  16
\epsilon^{1/2} \exp\{-x^2 /(4v_n^2)\} + \pr(|T_n-1|>\epsilon) .
$$
Moreover, for any $\delta >0$,
$$
\pr(|T_n -1|  > \epsilon)  \leq 4 \epsilon  ^{-1-\delta} \esp  \left[ |T_n
  -V_n^2|^{1+\delta} +|V_n^2 -1|^{1+\delta} \right] ,
$$
where $T_n -V_n^2$ is  a sum of Martingale differences. In the  same way as in
\citep{HallHeyde}, we obtain (as $\delta \leq 1$)
$$
\pr(|T_n -1| > \epsilon) \leq
C \epsilon ^{-1-\delta} \left[ \sum_{i=1} ^n \esp |Z_i|^{2+2\delta} +\esp |V_n^2 -1|^{1+\delta} \right] ,
$$
which concludes the proof.
\end{proof}

\end{document}